\documentclass[10pt]{amsart}
\usepackage[cp1251]{inputenc}
\usepackage[english]{babel}
\usepackage{amsmath}
\usepackage{amssymb}
\usepackage{amsfonts}

\usepackage[linktocpage=true, colorlinks=true, linkcolor=blue, citecolor=blue, urlcolor=blue]{hyperref}

\begin{document}
\renewcommand{\refname}{References}

\thispagestyle{empty}

\title[Expansions of Iterated Stratonovich Stochastic Integrals]
{Expansions of Iterated Stratonovich Stochastic Integrals 
from the Taylor--Stratonovich Expansion Based 
on Multiple Trigonometric Fourier Series. Comparison With the Milstein
Expansion}
\author[D.F. Kuznetsov]{Dmitriy F. Kuznetsov}
\address{Dmitriy Feliksovich Kuznetsov
\newline\hphantom{iii} Peter the Great Saint-Petersburg Polytechnic University,
\newline\hphantom{iii} Polytechnicheskaya ul., 29,
\newline\hphantom{iii} 195251, Saint-Petersburg, Russia}%
\email{sde\_kuznetsov@inbox.ru}

\maketitle {\small
\begin{quote}
\noindent{\sc Abstract.} 
The article is devoted to comparison of the Milstein 
expansion of iterated Stratonovich stochastic integrals with the method
of expansion of iterated stochastic integrals based on generalized
multiple Fourier series. We consider the practical material
connected with the expansions of iterated Stratonovich stochastic integrals
from the Taylor--Stratonovich expansion based on multiple trigonometric 
Fourier series.
The comparison of effectiveness of the Fourier--Legendre series as 
well as the trigonomertic
Fourier series for expansions of iterated Stratonovich stochastic integrals
is considered.
\medskip
\end{quote}
}

\vspace{12mm}


\setlength{\baselineskip}{1.5em}

\tableofcontents

\setlength{\baselineskip}{1.2em}


\vspace{5mm}

\section{Introduction}

\vspace{5mm}

Let $(\Omega,$ ${\rm F},$ ${\sf P})$ be a complete probability space, let 
$\{{\rm F}_t, t\in[0,T]\}$ be a nondecreasing right-continous 
family of $\sigma$-algebras of ${\rm F},$
and let ${\bf f}_t$ be a standard $m$-dimensional Wiener stochastic 
process, which is
${\rm F}_t$-measurable for any $t\in[0, T].$ We assume that the components
${\bf f}_{t}^{(i)}$ $(i=1,\ldots,m)$ of this process are independent. 
Consider
an Ito stochastic differential equation (SDE) in the integral form

\vspace{-1mm}
\begin{equation}
\label{1.5.2}
{\bf x}_t={\bf x}_0+\int\limits_0^t {\bf a}({\bf x}_{\tau},\tau)d\tau+
\int\limits_0^t B({\bf x}_{\tau},\tau)d{\bf f}_{\tau},\ \ \
{\bf x}_0={\bf x}(0,\omega).
\end{equation}

\vspace{3mm}
\noindent
Here ${\bf x}_t$ is some $n$-dimensional stochastic process 
satisfying the equation (\ref{1.5.2}). 
The nonrandom functions ${\bf a}:\ \mathbb{R}^n\times[0, T]\to\mathbb{R}^n$,
$B:\ \mathbb{R}^n\times[0, T]\to\mathbb{R}^{n\times m}$
guarantee the existence and uniqueness up to stochastic equivalence of 
a solution
of the equation (\ref{1.5.2}) \cite{1}. The second integral on the right-hand 
side of (\ref{1.5.2}) is 
interpreted as an Ito stochastic integral.
Let ${\bf x}_0$ be an $n$-dimensional random variable, which is 
${\rm F}_0$-measurable and 
${\sf M}\bigl\{\left|{\bf x}_0\right|^2\bigr\}<\infty$ 
(${\sf M}$ denotes a mathematical expectation).
We assume that
${\bf x}_0$ and ${\bf f}_t-{\bf f}_0$ are independent when $t>0.$

It is well known that one of the effective approaches 
to the numerical integration of 
Ito SDEs is an approach based on the Taylor--Ito and 
Taylor--Stratonovich expansions
\cite{Mi2}-\cite{Mi3}. The most important feature of such 
expansions is a presence in them of the so-called iterated
Ito and Stratonovich stochastic integrals, which play the key 
role for solving the 
problem of numerical integration of Ito SDEs and have the 
following form

\vspace{-1mm}
\begin{equation}
\label{ito}
J[\psi^{(k)}]_{T,t}=\int\limits_t^T\psi_k(t_k) \ldots \int\limits_t^{t_{2}}
\psi_1(t_1) d{\bf w}_{t_1}^{(i_1)}\ldots
d{\bf w}_{t_k}^{(i_k)},
\end{equation}

\vspace{1mm}
\begin{equation}
\label{str}
J^{*}[\psi^{(k)}]_{T,t}=
{\int\limits_t^{*}}^T
\psi_k(t_k)\ldots {\int\limits_t^{*}}^{t_2}
\psi_1(t_1) d{\bf w}_{t_1}^{(i_1)}\ldots d{\bf w}_{t_k}^{(i_k)},
\end{equation}

\vspace{3mm}
\noindent
where every $\psi_l(\tau)$ $(l=1,\ldots,k)$ is a
nonrandom function 
on $[t,T],$ ${\bf w}_{\tau}^{(i)}={\bf f}_{\tau}^{(i)}$
for $i=1,\ldots,m$ and
${\bf w}_{\tau}^{(0)}=\tau;$\
$i_1,\ldots,i_k = 0, 1,\ldots,m;$\ 

\vspace{-1mm}
$$
\int\limits\ \hbox{and}\ \int\limits^{*}
$$ 

\vspace{3mm}
\noindent
denote Ito and 
Stratonovich stochastic integrals,
respectively (in this paper, 
we use the definition of the Stratonovich stochastic integral from \cite{KlPl2}).

Note that $\psi_l(\tau)\equiv 1$ $(l=1,\ldots,k)$ and
$i_1,\ldots,i_k = 0, 1,\ldots,m$ for the classical 
Taylor--Ito and Taylor--Stratonovich expansions
\cite{Mi2}-\cite{Mi3} and 
$\psi_l(\tau)\equiv (t-\tau)^{q_l}$ ($l=1,\ldots,k$; 
$q_1,\ldots,q_k=0, 1, 2,\ldots $) and $i_1,\ldots,i_k = 1,\ldots,m$ 
for the unified Taylor--Ito and Taylor--Stratonovich expansions
\cite{k1}-\cite{2018axx}.

\vspace{5mm}

\section{Milstein Expansion and Method of 
Generatized Multiple Fourier Series}

\vspace{5mm}

Milstein G.N. proposed \cite{Mi2} (1988) an approach to the expansion 
of iterated stochastic integrals based on 
the trigonometric Fourier
expansion of the Brownian bridge 
process 
(version of the so-called 
Karhunen--Loeve expansion). 

Let us consider the Brownian bridge process \cite{Mi2}

\begin{equation}
\label{6.5.1}
{\bf f}_t-\frac{t}{\Delta}{\bf f}_{\Delta},\ \ \
t\in[0,\Delta],\ \ \ \Delta>0,
\end{equation}

\vspace{4mm}
\noindent
where ${\bf f}_t$ is a standard
$m$-dimensional Wiener process with independent components
${\bf f}^{(i)}_t$ $(i=1,\ldots,m).$

Consider the componentwise Karhunen--Loeve expansion 
of the process (\ref{6.5.1}) \cite{Mi2}

\begin{equation}
\label{6.5.2}
{\bf f}_t^{(i)}-\frac{t}{\Delta}{\bf f}_{\Delta}^{(i)}=
\frac{1}{2}a_{i,0}+\sum_{r=1}^{\infty}\left(
a_{i,r}{\rm cos}\frac{2\pi rt}{\Delta} +b_{i,r}{\rm sin}
\frac{2\pi rt}{\Delta}\right)
\end{equation}

\vspace{4mm}
\noindent
converging 
in the mean-square sense,
where

\vspace{-1mm}
$$
a_{i,r}=\frac{2}{\Delta} \int\limits_0^{\Delta}
\left({\bf f}_s^{(i)}-\frac{s}{\Delta}{\bf f}_{\Delta}^{(i)}\right)
{\rm cos}\frac{2\pi rs}{\Delta}ds,
$$

\vspace{1mm}
$$
b_{i,r}=\frac{2}{\Delta} \int\limits_0^{\Delta}
\left({\bf f}_s^{(i)}-\frac{s}{\Delta}{\bf f}_{\Delta}^{(i)}\right)
{\rm sin}\frac{2\pi rs}{\Delta}ds,
$$

\vspace{5mm}
\noindent
where $r=0, 1,\ldots;$ $i=1,\ldots,m.$ 

It is easy to demonstrate \cite{Mi2} that the random variables
$a_{i,r}, b_{i,r}$ 
are Gaussian ones and they satisfy the following relations

$$
{\sf M}\left\{a_{i,r}b_{i,r}\right\}=
{\sf M}\left\{a_{i,r}b_{i,k}\right\}=0,\ \ \
{\sf M}\left\{a_{i,r}a_{i,k}\right\}=
{\sf M}\left\{b_{i,r}b_{i,k}\right\}=0,
$$

\vspace{1mm}
$$
{\sf M}\left\{a_{i_1,r}a_{i_2,r}\right\}=
{\sf M}\left\{b_{i_1,r}b_{i_2,r}\right\}=0,\ \ \
{\sf M}\left\{a_{i,r}^2\right\}=
{\sf M}\left\{b_{i,r}^2\right\}=\frac{\Delta}{2\pi^2 r^2},
$$

\vspace{5mm}
\noindent
where $i, i_1, i_2=1,\ldots,m;$ $r\ne k;$ $i_1\ne i_2.$

According to (\ref{6.5.2}), we have

\begin{equation}
\label{6.5.7}
{\bf f}_t^{(i)}={\bf f}_{\Delta}^{(i)}\frac{t}{\Delta}+
\frac{1}{2}a_{i,0}+
\sum_{r=1}^{\infty}\left(
a_{i,r}{\rm cos}\frac{2\pi rt}{\Delta}+b_{i,r}{\rm sin}
\frac{2\pi rt}{\Delta}\right),
\end{equation}

\vspace{4mm}
\noindent
where the series
converges in the mean-square sense.

Note that the trigonometric functions are the 
eigenfunctions of the covariance operator
of the Brownian bridge process. That is why the basis functions are the
trigonometric functions in the considered approach.

In \cite{Mi2} 
Milstein G.N. proposed to expand (\ref{ito}) or (\ref{str})
(for the case $k=2$ and $\psi_1(s), \psi_2(s)\equiv 1$)
into iterated series of products
of standard Gaussian random variables by representing the Wiener
process as the series (\ref{6.5.7}).
To obtain the Milstein expansion of (\ref{ito}) or (\ref{str}), 
the truncated 
expansions (\ref{6.5.7}) of components of the Wiener 
process ${\bf f}_s$ must be
iteratively substituted in the single integrals, and the integrals
must be calculated, starting from the innermost integral.
This is a complicated procedure that obviously does not lead to a general
expansion of (\ref{ito}) or (\ref{str}) 
valid for an arbitrary multiplicity $k.$
For this reason, only expansions of simplest single, double, and triple
integrals (\ref{ito}), (\ref{str}) 
were obtained (see \cite{Mi2}-\cite{Mi3}).

At that, in \cite{Mi2}, \cite{Mi3} the case 
$\psi_1(s), \psi_2(s)\equiv 1$ and
$i_1, i_2=0, 1,\ldots,m$ is considered. In 
\cite{KlPl2}-\cite{Zapad-9} the attempt to consider the case 
$\psi_1(s), \psi_2(s), \psi_3(s)\equiv 1$ and 
$i_1, i_2, i_3=0, 1,\ldots,m$ is realized.

It should be noted that the authors of the works
\cite{KlPl2}
(Sect.~5.8, pp.~202--204), \cite{KPS} (pp.~82-84),
\cite{KPW} (pp.~438-439),  
\cite{Zapad-9} (pp.~263-264) use 
the Wong--Zakai approximation 
\cite{W-Z-1}-\cite{Watanabe} (without rigorous proof) within the frames
of the Milstein approach 
\cite{Mi2} based on the series expansion 
of the Brownian bridge process. See discussion in Sect.~7 of 
this paper for details.

Let us consider an another approach to the expansion of iterated
stochastic integrals \cite{2006}-\cite{arxiv-8}, 
which is reffered to as the method of generalized
multiple Fourier series.

Suppose that every $\psi_l(\tau)$ $(l=1,\ldots,k)$ is a 
nonrandom function from the space $L_2([t, T])$. 
Define the following function on the hypercube $[t, T]^k$

\begin{equation}
\label{ppp}
K(t_1,\ldots,t_k)=
\begin{cases}
\psi_1(t_1)\ldots \psi_k(t_k),\ t_1<\ldots<t_k\\
~\\
~\\
0,\ \hbox{\rm otherwise}
\end{cases}
=\ \
\prod\limits_{l=1}^k
\psi_l(t_l)\ \prod\limits_{l=1}^{k-1}{\bf 1}_{\{t_l<t_{l+1}\}},\ 
\end{equation}

\vspace{5mm}
\noindent
where $t_1,\ldots,t_k\in [t, T]$ $(k\ge 2)$ and 
$K(t_1)\equiv\psi_1(t_1)$ for $t_1\in[t, T].$ Here 
${\bf 1}_A$ denotes the indicator of the set $A$.

Suppose that $\{\phi_j(x)\}_{j=0}^{\infty}$
is a complete orthonormal system of functions in 
the space $L_2([t, T])$.

The function $K(t_1,\ldots,t_k)$ belongs to the space $L_2([t, T]^k).$
At this situation it is well known that the generalized 
multiple Fourier series 
of $K(t_1,\ldots,t_k)\in L_2([t, T]^k)$ is converging 
to $K(t_1,\ldots,t_k)$ in the hypercube $[t, T]^k$ in 
the mean-square sense, i.e.

\begin{equation}
\label{sos1z}
\hbox{\vtop{\offinterlineskip\halign{
\hfil#\hfil\cr
{\rm lim}\cr
$\stackrel{}{{}_{p_1,\ldots,p_k\to \infty}}$\cr
}} }\Biggl\Vert
K(t_1,\ldots,t_k)-
\sum_{j_1=0}^{p_1}\ldots \sum_{j_k=0}^{p_k}
C_{j_k\ldots j_1}\prod_{l=1}^{k} \phi_{j_l}(t_l)\Biggr
\Vert_{L_2([t,T]^k)}=0,
\end{equation}

\vspace{3mm}
\noindent
where
\begin{equation}
\label{ppppa}
C_{j_k\ldots j_1}=\int\limits_{[t,T]^k}
K(t_1,\ldots,t_k)\prod_{l=1}^{k}\phi_{j_l}(t_l)dt_1\ldots dt_k
\end{equation}

\vspace{4mm}
\noindent
is the Fourier coefficient and

\vspace{-2mm}
$$
\left\Vert f\right\Vert_{L_2([t,T]^k)}=\left(\int\limits_{[t,T]^k}
f^2(t_1,\ldots,t_k)dt_1\ldots dt_k\right)^{1/2}.
$$

\vspace{4mm}

Consider the partition $\{\tau_j\}_{j=0}^N$ of $[t, T]$ such that

\begin{equation}
\label{1111}
t=\tau_0<\ldots <\tau_N=T,\ \ \
\Delta_N=
\hbox{\vtop{\offinterlineskip\halign{
\hfil#\hfil\cr
{\rm max}\cr
$\stackrel{}{{}_{0\le j\le N-1}}$\cr
}} }\Delta\tau_j\to 0\ \ \hbox{if}\ \ N\to \infty,\ \ \ 
\Delta\tau_j=\tau_{j+1}-\tau_j.
\end{equation}

\vspace{5mm}

{\bf Theorem 1} \cite{2006} (2006), \cite{2007-1}-\cite{arxiv-8}. 
{\it Suppose that
every $\psi_l(\tau)$ $(l=1,\ldots, k)$ is a continuous non-random
function on 
$[t, T]$ and
$\{\phi_j(x)\}_{j=0}^{\infty}$ is a complete orthonormal system  
of continuous functions in the space $L_2([t,T]).$ 
Then

\vspace{1mm}
$$
J[\psi^{(k)}]_{T,t}\  =\ 
\hbox{\vtop{\offinterlineskip\halign{
\hfil#\hfil\cr
{\rm l.i.m.}\cr
$\stackrel{}{{}_{p_1,\ldots,p_k\to \infty}}$\cr
}} }\sum_{j_1=0}^{p_1}\ldots\sum_{j_k=0}^{p_k}
C_{j_k\ldots j_1}\Biggl(
\prod_{l=1}^k\zeta_{j_l}^{(i_l)}\ -
\Biggr.
$$

\vspace{3mm}
\begin{equation}
\label{tyyy}
-\ \Biggl.
\hbox{\vtop{\offinterlineskip\halign{
\hfil#\hfil\cr
{\rm l.i.m.}\cr
$\stackrel{}{{}_{N\to \infty}}$\cr
}} }\sum_{(l_1,\ldots,l_k)\in {\rm G}_k}
\phi_{j_{1}}(\tau_{l_1})
\Delta{\bf w}_{\tau_{l_1}}^{(i_1)}\ldots
\phi_{j_{k}}(\tau_{l_k})
\Delta{\bf w}_{\tau_{l_k}}^{(i_k)}\Biggr),
\end{equation}

\vspace{6mm}
\noindent
where

\vspace{-2mm}
$$
{\rm G}_k={\rm H}_k\backslash{\rm L}_k,\ \ \
{\rm H}_k=\{(l_1,\ldots,l_k):\ l_1,\ldots,l_k=0,\ 1,\ldots,N-1\},
$$

\vspace{1mm}
$$
{\rm L}_k=\{(l_1,\ldots,l_k):\ l_1,\ldots,l_k=0,\ 1,\ldots,N-1;\
l_g\ne l_r\ (g\ne r);\ g, r=1,\ldots,k\},
$$

\vspace{6mm}
\noindent
${\rm l.i.m.}$ is a limit in the mean-square sense,
$i_1,\ldots,i_k=0,1,\ldots,m,$

\vspace{-1mm}
\begin{equation}
\label{rr23}
\zeta_{j}^{(i)}=
\int\limits_t^T \phi_{j}(s) d{\bf w}_s^{(i)}
\end{equation} 

\vspace{3mm}
\noindent
are independent standard Gaussian random variables
for various
$i$ or $j$ {\rm(}if $i\ne 0${\rm),}
$C_{j_k\ldots j_1}$ is the Fourier coefficient {\rm(\ref{ppppa}),}
$\Delta{\bf w}_{\tau_{j}}^{(i)}=
{\bf w}_{\tau_{j+1}}^{(i)}-{\bf w}_{\tau_{j}}^{(i)}$
$(i=0, 1,\ldots,m),$\
$\left\{\tau_{j}\right\}_{j=0}^{N}$ is the partition of
$[t,T],$ which satisfies the condition {\rm (\ref{1111})}.}

\vspace{2mm}

In order to evaluate the significance of Theorem 1 for practice we will
demonstrate its transformed particular cases for 
$k=1,\ldots,6$ \cite{2006}-\cite{arxiv-8}

\vspace{1mm}
\begin{equation}
\label{a1}
J[\psi^{(1)}]_{T,t}
=\hbox{\vtop{\offinterlineskip\halign{
\hfil#\hfil\cr
{\rm l.i.m.}\cr
$\stackrel{}{{}_{p_1\to \infty}}$\cr
}} }\sum_{j_1=0}^{p_1}
C_{j_1}\zeta_{j_1}^{(i_1)},
\end{equation}

\vspace{4mm}
\begin{equation}
\label{a2}
J[\psi^{(2)}]_{T,t}
=\hbox{\vtop{\offinterlineskip\halign{
\hfil#\hfil\cr
{\rm l.i.m.}\cr
$\stackrel{}{{}_{p_1,p_2\to \infty}}$\cr
}} }\sum_{j_1=0}^{p_1}\sum_{j_2=0}^{p_2}
C_{j_2j_1}\Biggl(\zeta_{j_1}^{(i_1)}\zeta_{j_2}^{(i_2)}
-{\bf 1}_{\{i_1=i_2\ne 0\}}
{\bf 1}_{\{j_1=j_2\}}\Biggr),
\end{equation}

\vspace{6mm}
$$
J[\psi^{(3)}]_{T,t}=
\hbox{\vtop{\offinterlineskip\halign{
\hfil#\hfil\cr
{\rm l.i.m.}\cr
$\stackrel{}{{}_{p_1,\ldots,p_3\to \infty}}$\cr
}} }\sum_{j_1=0}^{p_1}\sum_{j_2=0}^{p_2}\sum_{j_3=0}^{p_3}
C_{j_3j_2j_1}\Biggl(
\zeta_{j_1}^{(i_1)}\zeta_{j_2}^{(i_2)}\zeta_{j_3}^{(i_3)}
-\Biggr.
$$

\vspace{1mm}
\begin{equation}
\label{a3}
\Biggl.-{\bf 1}_{\{i_1=i_2\ne 0\}}
{\bf 1}_{\{j_1=j_2\}}
\zeta_{j_3}^{(i_3)}
-{\bf 1}_{\{i_2=i_3\ne 0\}}
{\bf 1}_{\{j_2=j_3\}}
\zeta_{j_1}^{(i_1)}-
{\bf 1}_{\{i_1=i_3\ne 0\}}
{\bf 1}_{\{j_1=j_3\}}
\zeta_{j_2}^{(i_2)}\Biggr),
\end{equation}

\vspace{6mm}

$$
J[\psi^{(4)}]_{T,t}
=
\hbox{\vtop{\offinterlineskip\halign{
\hfil#\hfil\cr
{\rm l.i.m.}\cr
$\stackrel{}{{}_{p_1,\ldots,p_4\to \infty}}$\cr
}} }\sum_{j_1=0}^{p_1}\ldots\sum_{j_4=0}^{p_4}
C_{j_4\ldots j_1}\Biggl(
\prod_{l=1}^4\zeta_{j_l}^{(i_l)}
\Biggr.
-
$$
$$
-
{\bf 1}_{\{i_1=i_2\ne 0\}}
{\bf 1}_{\{j_1=j_2\}}
\zeta_{j_3}^{(i_3)}
\zeta_{j_4}^{(i_4)}
-
{\bf 1}_{\{i_1=i_3\ne 0\}}
{\bf 1}_{\{j_1=j_3\}}
\zeta_{j_2}^{(i_2)}
\zeta_{j_4}^{(i_4)}-
$$
$$
-
{\bf 1}_{\{i_1=i_4\ne 0\}}
{\bf 1}_{\{j_1=j_4\}}
\zeta_{j_2}^{(i_2)}
\zeta_{j_3}^{(i_3)}
-
{\bf 1}_{\{i_2=i_3\ne 0\}}
{\bf 1}_{\{j_2=j_3\}}
\zeta_{j_1}^{(i_1)}
\zeta_{j_4}^{(i_4)}-
$$
$$
-
{\bf 1}_{\{i_2=i_4\ne 0\}}
{\bf 1}_{\{j_2=j_4\}}
\zeta_{j_1}^{(i_1)}
\zeta_{j_3}^{(i_3)}
-
{\bf 1}_{\{i_3=i_4\ne 0\}}
{\bf 1}_{\{j_3=j_4\}}
\zeta_{j_1}^{(i_1)}
\zeta_{j_2}^{(i_2)}+
$$
$$
+
{\bf 1}_{\{i_1=i_2\ne 0\}}
{\bf 1}_{\{j_1=j_2\}}
{\bf 1}_{\{i_3=i_4\ne 0\}}
{\bf 1}_{\{j_3=j_4\}}
+
$$
$$
+
{\bf 1}_{\{i_1=i_3\ne 0\}}
{\bf 1}_{\{j_1=j_3\}}
{\bf 1}_{\{i_2=i_4\ne 0\}}
{\bf 1}_{\{j_2=j_4\}}+
$$
\begin{equation}
\label{a4}
+\Biggl.
{\bf 1}_{\{i_1=i_4\ne 0\}}
{\bf 1}_{\{j_1=j_4\}}
{\bf 1}_{\{i_2=i_3\ne 0\}}
{\bf 1}_{\{j_2=j_3\}}\Biggr),
\end{equation}

\vspace{8mm}

$$
J[\psi^{(5)}]_{T,t}
=\hbox{\vtop{\offinterlineskip\halign{
\hfil#\hfil\cr
{\rm l.i.m.}\cr
$\stackrel{}{{}_{p_1,\ldots,p_5\to \infty}}$\cr
}} }\sum_{j_1=0}^{p_1}\ldots\sum_{j_5=0}^{p_5}
C_{j_5\ldots j_1}\Biggl(
\prod_{l=1}^5\zeta_{j_l}^{(i_l)}
-\Biggr.
$$
$$
-
{\bf 1}_{\{i_1=i_2\ne 0\}}
{\bf 1}_{\{j_1=j_2\}}
\zeta_{j_3}^{(i_3)}
\zeta_{j_4}^{(i_4)}
\zeta_{j_5}^{(i_5)}-
{\bf 1}_{\{i_1=i_3\ne 0\}}
{\bf 1}_{\{j_1=j_3\}}
\zeta_{j_2}^{(i_2)}
\zeta_{j_4}^{(i_4)}
\zeta_{j_5}^{(i_5)}-
$$
$$
-
{\bf 1}_{\{i_1=i_4\ne 0\}}
{\bf 1}_{\{j_1=j_4\}}
\zeta_{j_2}^{(i_2)}
\zeta_{j_3}^{(i_3)}
\zeta_{j_5}^{(i_5)}-
{\bf 1}_{\{i_1=i_5\ne 0\}}
{\bf 1}_{\{j_1=j_5\}}
\zeta_{j_2}^{(i_2)}
\zeta_{j_3}^{(i_3)}
\zeta_{j_4}^{(i_4)}-
$$
$$
-
{\bf 1}_{\{i_2=i_3\ne 0\}}
{\bf 1}_{\{j_2=j_3\}}
\zeta_{j_1}^{(i_1)}
\zeta_{j_4}^{(i_4)}
\zeta_{j_5}^{(i_5)}-
{\bf 1}_{\{i_2=i_4\ne 0\}}
{\bf 1}_{\{j_2=j_4\}}
\zeta_{j_1}^{(i_1)}
\zeta_{j_3}^{(i_3)}
\zeta_{j_5}^{(i_5)}-
$$
$$
-
{\bf 1}_{\{i_2=i_5\ne 0\}}
{\bf 1}_{\{j_2=j_5\}}
\zeta_{j_1}^{(i_1)}
\zeta_{j_3}^{(i_3)}
\zeta_{j_4}^{(i_4)}
-{\bf 1}_{\{i_3=i_4\ne 0\}}
{\bf 1}_{\{j_3=j_4\}}
\zeta_{j_1}^{(i_1)}
\zeta_{j_2}^{(i_2)}
\zeta_{j_5}^{(i_5)}-
$$
$$
-
{\bf 1}_{\{i_3=i_5\ne 0\}}
{\bf 1}_{\{j_3=j_5\}}
\zeta_{j_1}^{(i_1)}
\zeta_{j_2}^{(i_2)}
\zeta_{j_4}^{(i_4)}
-{\bf 1}_{\{i_4=i_5\ne 0\}}
{\bf 1}_{\{j_4=j_5\}}
\zeta_{j_1}^{(i_1)}
\zeta_{j_2}^{(i_2)}
\zeta_{j_3}^{(i_3)}+
$$
$$
+
{\bf 1}_{\{i_1=i_2\ne 0\}}
{\bf 1}_{\{j_1=j_2\}}
{\bf 1}_{\{i_3=i_4\ne 0\}}
{\bf 1}_{\{j_3=j_4\}}\zeta_{j_5}^{(i_5)}+
{\bf 1}_{\{i_1=i_2\ne 0\}}
{\bf 1}_{\{j_1=j_2\}}
{\bf 1}_{\{i_3=i_5\ne 0\}}
{\bf 1}_{\{j_3=j_5\}}\zeta_{j_4}^{(i_4)}+
$$
$$
+
{\bf 1}_{\{i_1=i_2\ne 0\}}
{\bf 1}_{\{j_1=j_2\}}
{\bf 1}_{\{i_4=i_5\ne 0\}}
{\bf 1}_{\{j_4=j_5\}}\zeta_{j_3}^{(i_3)}+
{\bf 1}_{\{i_1=i_3\ne 0\}}
{\bf 1}_{\{j_1=j_3\}}
{\bf 1}_{\{i_2=i_4\ne 0\}}
{\bf 1}_{\{j_2=j_4\}}\zeta_{j_5}^{(i_5)}+
$$
$$
+
{\bf 1}_{\{i_1=i_3\ne 0\}}
{\bf 1}_{\{j_1=j_3\}}
{\bf 1}_{\{i_2=i_5\ne 0\}}
{\bf 1}_{\{j_2=j_5\}}\zeta_{j_4}^{(i_4)}+
{\bf 1}_{\{i_1=i_3\ne 0\}}
{\bf 1}_{\{j_1=j_3\}}
{\bf 1}_{\{i_4=i_5\ne 0\}}
{\bf 1}_{\{j_4=j_5\}}\zeta_{j_2}^{(i_2)}+
$$
$$
+
{\bf 1}_{\{i_1=i_4\ne 0\}}
{\bf 1}_{\{j_1=j_4\}}
{\bf 1}_{\{i_2=i_3\ne 0\}}
{\bf 1}_{\{j_2=j_3\}}\zeta_{j_5}^{(i_5)}+
{\bf 1}_{\{i_1=i_4\ne 0\}}
{\bf 1}_{\{j_1=j_4\}}
{\bf 1}_{\{i_2=i_5\ne 0\}}
{\bf 1}_{\{j_2=j_5\}}\zeta_{j_3}^{(i_3)}+
$$
$$
+
{\bf 1}_{\{i_1=i_4\ne 0\}}
{\bf 1}_{\{j_1=j_4\}}
{\bf 1}_{\{i_3=i_5\ne 0\}}
{\bf 1}_{\{j_3=j_5\}}\zeta_{j_2}^{(i_2)}+
{\bf 1}_{\{i_1=i_5\ne 0\}}
{\bf 1}_{\{j_1=j_5\}}
{\bf 1}_{\{i_2=i_3\ne 0\}}
{\bf 1}_{\{j_2=j_3\}}\zeta_{j_4}^{(i_4)}+
$$
$$
+
{\bf 1}_{\{i_1=i_5\ne 0\}}
{\bf 1}_{\{j_1=j_5\}}
{\bf 1}_{\{i_2=i_4\ne 0\}}
{\bf 1}_{\{j_2=j_4\}}\zeta_{j_3}^{(i_3)}+
{\bf 1}_{\{i_1=i_5\ne 0\}}
{\bf 1}_{\{j_1=j_5\}}
{\bf 1}_{\{i_3=i_4\ne 0\}}
{\bf 1}_{\{j_3=j_4\}}\zeta_{j_2}^{(i_2)}+
$$
$$
+
{\bf 1}_{\{i_2=i_3\ne 0\}}
{\bf 1}_{\{j_2=j_3\}}
{\bf 1}_{\{i_4=i_5\ne 0\}}
{\bf 1}_{\{j_4=j_5\}}\zeta_{j_1}^{(i_1)}+
{\bf 1}_{\{i_2=i_4\ne 0\}}
{\bf 1}_{\{j_2=j_4\}}
{\bf 1}_{\{i_3=i_5\ne 0\}}
{\bf 1}_{\{j_3=j_5\}}\zeta_{j_1}^{(i_1)}+
$$
\begin{equation}
\label{a5}
+\Biggl.
{\bf 1}_{\{i_2=i_5\ne 0\}}
{\bf 1}_{\{j_2=j_5\}}
{\bf 1}_{\{i_3=i_4\ne 0\}}
{\bf 1}_{\{j_3=j_4\}}\zeta_{j_1}^{(i_1)}\Biggr),
\end{equation}

\vspace{9mm}

$$
J[\psi^{(6)}]_{T,t}
=\hbox{\vtop{\offinterlineskip\halign{
\hfil#\hfil\cr
{\rm l.i.m.}\cr
$\stackrel{}{{}_{p_1,\ldots,p_6\to \infty}}$\cr
}} }\sum_{j_1=0}^{p_1}\ldots\sum_{j_6=0}^{p_6}
C_{j_6\ldots j_1}\Biggl(
\prod_{l=1}^6
\zeta_{j_l}^{(i_l)}
-\Biggr.
$$
$$
-
{\bf 1}_{\{i_1=i_6\ne 0\}}
{\bf 1}_{\{j_1=j_6\}}
\zeta_{j_2}^{(i_2)}
\zeta_{j_3}^{(i_3)}
\zeta_{j_4}^{(i_4)}
\zeta_{j_5}^{(i_5)}-
{\bf 1}_{\{i_2=i_6\ne 0\}}
{\bf 1}_{\{j_2=j_6\}}
\zeta_{j_1}^{(i_1)}
\zeta_{j_3}^{(i_3)}
\zeta_{j_4}^{(i_4)}
\zeta_{j_5}^{(i_5)}-
$$
$$
-
{\bf 1}_{\{i_3=i_6\ne 0\}}
{\bf 1}_{\{j_3=j_6\}}
\zeta_{j_1}^{(i_1)}
\zeta_{j_2}^{(i_2)}
\zeta_{j_4}^{(i_4)}
\zeta_{j_5}^{(i_5)}-
{\bf 1}_{\{i_4=i_6\ne 0\}}
{\bf 1}_{\{j_4=j_6\}}
\zeta_{j_1}^{(i_1)}
\zeta_{j_2}^{(i_2)}
\zeta_{j_3}^{(i_3)}
\zeta_{j_5}^{(i_5)}-
$$
$$
-
{\bf 1}_{\{i_5=i_6\ne 0\}}
{\bf 1}_{\{j_5=j_6\}}
\zeta_{j_1}^{(i_1)}
\zeta_{j_2}^{(i_2)}
\zeta_{j_3}^{(i_3)}
\zeta_{j_4}^{(i_4)}-
{\bf 1}_{\{i_1=i_2\ne 0\}}
{\bf 1}_{\{j_1=j_2\}}
\zeta_{j_3}^{(i_3)}
\zeta_{j_4}^{(i_4)}
\zeta_{j_5}^{(i_5)}
\zeta_{j_6}^{(i_6)}-
$$
$$
-
{\bf 1}_{\{i_1=i_3\ne 0\}}
{\bf 1}_{\{j_1=j_3\}}
\zeta_{j_2}^{(i_2)}
\zeta_{j_4}^{(i_4)}
\zeta_{j_5}^{(i_5)}
\zeta_{j_6}^{(i_6)}-
{\bf 1}_{\{i_1=i_4\ne 0\}}
{\bf 1}_{\{j_1=j_4\}}
\zeta_{j_2}^{(i_2)}
\zeta_{j_3}^{(i_3)}
\zeta_{j_5}^{(i_5)}
\zeta_{j_6}^{(i_6)}-
$$
$$
-
{\bf 1}_{\{i_1=i_5\ne 0\}}
{\bf 1}_{\{j_1=j_5\}}
\zeta_{j_2}^{(i_2)}
\zeta_{j_3}^{(i_3)}
\zeta_{j_4}^{(i_4)}
\zeta_{j_6}^{(i_6)}-
{\bf 1}_{\{i_2=i_3\ne 0\}}
{\bf 1}_{\{j_2=j_3\}}
\zeta_{j_1}^{(i_1)}
\zeta_{j_4}^{(i_4)}
\zeta_{j_5}^{(i_5)}
\zeta_{j_6}^{(i_6)}-
$$
$$
-
{\bf 1}_{\{i_2=i_4\ne 0\}}
{\bf 1}_{\{j_2=j_4\}}
\zeta_{j_1}^{(i_1)}
\zeta_{j_3}^{(i_3)}
\zeta_{j_5}^{(i_5)}
\zeta_{j_6}^{(i_6)}-
{\bf 1}_{\{i_2=i_5\ne 0\}}
{\bf 1}_{\{j_2=j_5\}}
\zeta_{j_1}^{(i_1)}
\zeta_{j_3}^{(i_3)}
\zeta_{j_4}^{(i_4)}
\zeta_{j_6}^{(i_6)}-
$$
$$
-
{\bf 1}_{\{i_3=i_4\ne 0\}}
{\bf 1}_{\{j_3=j_4\}}
\zeta_{j_1}^{(i_1)}
\zeta_{j_2}^{(i_2)}
\zeta_{j_5}^{(i_5)}
\zeta_{j_6}^{(i_6)}-
{\bf 1}_{\{i_3=i_5\ne 0\}}
{\bf 1}_{\{j_3=j_5\}}
\zeta_{j_1}^{(i_1)}
\zeta_{j_2}^{(i_2)}
\zeta_{j_4}^{(i_4)}
\zeta_{j_6}^{(i_6)}-
$$
$$
-
{\bf 1}_{\{i_4=i_5\ne 0\}}
{\bf 1}_{\{j_4=j_5\}}
\zeta_{j_1}^{(i_1)}
\zeta_{j_2}^{(i_2)}
\zeta_{j_3}^{(i_3)}
\zeta_{j_6}^{(i_6)}+
$$
$$
+
{\bf 1}_{\{i_1=i_2\ne 0\}}
{\bf 1}_{\{j_1=j_2\}}
{\bf 1}_{\{i_3=i_4\ne 0\}}
{\bf 1}_{\{j_3=j_4\}}
\zeta_{j_5}^{(i_5)}
\zeta_{j_6}^{(i_6)}+
{\bf 1}_{\{i_1=i_2\ne 0\}}
{\bf 1}_{\{j_1=j_2\}}
{\bf 1}_{\{i_3=i_5\ne 0\}}
{\bf 1}_{\{j_3=j_5\}}
\zeta_{j_4}^{(i_4)}
\zeta_{j_6}^{(i_6)}+
$$
$$
+
{\bf 1}_{\{i_1=i_2\ne 0\}}
{\bf 1}_{\{j_1=j_2\}}
{\bf 1}_{\{i_4=i_5\ne 0\}}
{\bf 1}_{\{j_4=j_5\}}
\zeta_{j_3}^{(i_3)}
\zeta_{j_6}^{(i_6)}
+
{\bf 1}_{\{i_1=i_3\ne 0\}}
{\bf 1}_{\{j_1=j_3\}}
{\bf 1}_{\{i_2=i_4\ne 0\}}
{\bf 1}_{\{j_2=j_4\}}
\zeta_{j_5}^{(i_5)}
\zeta_{j_6}^{(i_6)}+
$$
$$
+
{\bf 1}_{\{i_1=i_3\ne 0\}}
{\bf 1}_{\{j_1=j_3\}}
{\bf 1}_{\{i_2=i_5\ne 0\}}
{\bf 1}_{\{j_2=j_5\}}
\zeta_{j_4}^{(i_4)}
\zeta_{j_6}^{(i_6)}
+{\bf 1}_{\{i_1=i_3\ne 0\}}
{\bf 1}_{\{j_1=j_3\}}
{\bf 1}_{\{i_4=i_5\ne 0\}}
{\bf 1}_{\{j_4=j_5\}}
\zeta_{j_2}^{(i_2)}
\zeta_{j_6}^{(i_6)}+
$$
$$
+
{\bf 1}_{\{i_1=i_4\ne 0\}}
{\bf 1}_{\{j_1=j_4\}}
{\bf 1}_{\{i_2=i_3\ne 0\}}
{\bf 1}_{\{j_2=j_3\}}
\zeta_{j_5}^{(i_5)}
\zeta_{j_6}^{(i_6)}
+
{\bf 1}_{\{i_1=i_4\ne 0\}}
{\bf 1}_{\{j_1=j_4\}}
{\bf 1}_{\{i_2=i_5\ne 0\}}
{\bf 1}_{\{j_2=j_5\}}
\zeta_{j_3}^{(i_3)}
\zeta_{j_6}^{(i_6)}+
$$
$$
+
{\bf 1}_{\{i_1=i_4\ne 0\}}
{\bf 1}_{\{j_1=j_4\}}
{\bf 1}_{\{i_3=i_5\ne 0\}}
{\bf 1}_{\{j_3=j_5\}}
\zeta_{j_2}^{(i_2)}
\zeta_{j_6}^{(i_6)}
+
{\bf 1}_{\{i_1=i_5\ne 0\}}
{\bf 1}_{\{j_1=j_5\}}
{\bf 1}_{\{i_2=i_3\ne 0\}}
{\bf 1}_{\{j_2=j_3\}}
\zeta_{j_4}^{(i_4)}
\zeta_{j_6}^{(i_6)}+
$$
$$
+
{\bf 1}_{\{i_1=i_5\ne 0\}}
{\bf 1}_{\{j_1=j_5\}}
{\bf 1}_{\{i_2=i_4\ne 0\}}
{\bf 1}_{\{j_2=j_4\}}
\zeta_{j_3}^{(i_3)}
\zeta_{j_6}^{(i_6)}
+
{\bf 1}_{\{i_1=i_5\ne 0\}}
{\bf 1}_{\{j_1=j_5\}}
{\bf 1}_{\{i_3=i_4\ne 0\}}
{\bf 1}_{\{j_3=j_4\}}
\zeta_{j_2}^{(i_2)}
\zeta_{j_6}^{(i_6)}+
$$
$$
+
{\bf 1}_{\{i_2=i_3\ne 0\}}
{\bf 1}_{\{j_2=j_3\}}
{\bf 1}_{\{i_4=i_5\ne 0\}}
{\bf 1}_{\{j_4=j_5\}}
\zeta_{j_1}^{(i_1)}
\zeta_{j_6}^{(i_6)}
+
{\bf 1}_{\{i_2=i_4\ne 0\}}
{\bf 1}_{\{j_2=j_4\}}
{\bf 1}_{\{i_3=i_5\ne 0\}}
{\bf 1}_{\{j_3=j_5\}}
\zeta_{j_1}^{(i_1)}
\zeta_{j_6}^{(i_6)}+
$$
$$
+
{\bf 1}_{\{i_2=i_5\ne 0\}}
{\bf 1}_{\{j_2=j_5\}}
{\bf 1}_{\{i_3=i_4\ne 0\}}
{\bf 1}_{\{j_3=j_4\}}
\zeta_{j_1}^{(i_1)}
\zeta_{j_6}^{(i_6)}
+
{\bf 1}_{\{i_6=i_1\ne 0\}}
{\bf 1}_{\{j_6=j_1\}}
{\bf 1}_{\{i_3=i_4\ne 0\}}
{\bf 1}_{\{j_3=j_4\}}
\zeta_{j_2}^{(i_2)}
\zeta_{j_5}^{(i_5)}+
$$
$$
+
{\bf 1}_{\{i_6=i_1\ne 0\}}
{\bf 1}_{\{j_6=j_1\}}
{\bf 1}_{\{i_3=i_5\ne 0\}}
{\bf 1}_{\{j_3=j_5\}}
\zeta_{j_2}^{(i_2)}
\zeta_{j_4}^{(i_4)}
+
{\bf 1}_{\{i_6=i_1\ne 0\}}
{\bf 1}_{\{j_6=j_1\}}
{\bf 1}_{\{i_2=i_5\ne 0\}}
{\bf 1}_{\{j_2=j_5\}}
\zeta_{j_3}^{(i_3)}
\zeta_{j_4}^{(i_4)}+
$$
$$
+
{\bf 1}_{\{i_6=i_1\ne 0\}}
{\bf 1}_{\{j_6=j_1\}}
{\bf 1}_{\{i_2=i_4\ne 0\}}
{\bf 1}_{\{j_2=j_4\}}
\zeta_{j_3}^{(i_3)}
\zeta_{j_5}^{(i_5)}
+
{\bf 1}_{\{i_6=i_1\ne 0\}}
{\bf 1}_{\{j_6=j_1\}}
{\bf 1}_{\{i_4=i_5\ne 0\}}
{\bf 1}_{\{j_4=j_5\}}
\zeta_{j_2}^{(i_2)}
\zeta_{j_3}^{(i_3)}+
$$
$$
+
{\bf 1}_{\{i_6=i_1\ne 0\}}
{\bf 1}_{\{j_6=j_1\}}
{\bf 1}_{\{i_2=i_3\ne 0\}}
{\bf 1}_{\{j_2=j_3\}}
\zeta_{j_4}^{(i_4)}
\zeta_{j_5}^{(i_5)}
+
{\bf 1}_{\{i_6=i_2\ne 0\}}
{\bf 1}_{\{j_6=j_2\}}
{\bf 1}_{\{i_3=i_5\ne 0\}}
{\bf 1}_{\{j_3=j_5\}}
\zeta_{j_1}^{(i_1)}
\zeta_{j_4}^{(i_4)}+
$$
$$
+
{\bf 1}_{\{i_6=i_2\ne 0\}}
{\bf 1}_{\{j_6=j_2\}}
{\bf 1}_{\{i_4=i_5\ne 0\}}
{\bf 1}_{\{j_4=j_5\}}
\zeta_{j_1}^{(i_1)}
\zeta_{j_3}^{(i_3)}
+
{\bf 1}_{\{i_6=i_2\ne 0\}}
{\bf 1}_{\{j_6=j_2\}}
{\bf 1}_{\{i_3=i_4\ne 0\}}
{\bf 1}_{\{j_3=j_4\}}
\zeta_{j_1}^{(i_1)}
\zeta_{j_5}^{(i_5)}+
$$
$$
+
{\bf 1}_{\{i_6=i_2\ne 0\}}
{\bf 1}_{\{j_6=j_2\}}
{\bf 1}_{\{i_1=i_5\ne 0\}}
{\bf 1}_{\{j_1=j_5\}}
\zeta_{j_3}^{(i_3)}
\zeta_{j_4}^{(i_4)}
+
{\bf 1}_{\{i_6=i_2\ne 0\}}
{\bf 1}_{\{j_6=j_2\}}
{\bf 1}_{\{i_1=i_4\ne 0\}}
{\bf 1}_{\{j_1=j_4\}}
\zeta_{j_3}^{(i_3)}
\zeta_{j_5}^{(i_5)}+
$$
$$
+
{\bf 1}_{\{i_6=i_2\ne 0\}}
{\bf 1}_{\{j_6=j_2\}}
{\bf 1}_{\{i_1=i_3\ne 0\}}
{\bf 1}_{\{j_1=j_3\}}
\zeta_{j_4}^{(i_4)}
\zeta_{j_5}^{(i_5)}
+
{\bf 1}_{\{i_6=i_3\ne 0\}}
{\bf 1}_{\{j_6=j_3\}}
{\bf 1}_{\{i_2=i_5\ne 0\}}
{\bf 1}_{\{j_2=j_5\}}
\zeta_{j_1}^{(i_1)}
\zeta_{j_4}^{(i_4)}+
$$
$$
+
{\bf 1}_{\{i_6=i_3\ne 0\}}
{\bf 1}_{\{j_6=j_3\}}
{\bf 1}_{\{i_4=i_5\ne 0\}}
{\bf 1}_{\{j_4=j_5\}}
\zeta_{j_1}^{(i_1)}
\zeta_{j_2}^{(i_2)}
+
{\bf 1}_{\{i_6=i_3\ne 0\}}
{\bf 1}_{\{j_6=j_3\}}
{\bf 1}_{\{i_2=i_4\ne 0\}}
{\bf 1}_{\{j_2=j_4\}}
\zeta_{j_1}^{(i_1)}
\zeta_{j_5}^{(i_5)}+
$$
$$
+
{\bf 1}_{\{i_6=i_3\ne 0\}}
{\bf 1}_{\{j_6=j_3\}}
{\bf 1}_{\{i_1=i_5\ne 0\}}
{\bf 1}_{\{j_1=j_5\}}
\zeta_{j_2}^{(i_2)}
\zeta_{j_4}^{(i_4)}
+
{\bf 1}_{\{i_6=i_3\ne 0\}}
{\bf 1}_{\{j_6=j_3\}}
{\bf 1}_{\{i_1=i_4\ne 0\}}
{\bf 1}_{\{j_1=j_4\}}
\zeta_{j_2}^{(i_2)}
\zeta_{j_5}^{(i_5)}+
$$
$$
+
{\bf 1}_{\{i_6=i_3\ne 0\}}
{\bf 1}_{\{j_6=j_3\}}
{\bf 1}_{\{i_1=i_2\ne 0\}}
{\bf 1}_{\{j_1=j_2\}}
\zeta_{j_4}^{(i_4)}
\zeta_{j_5}^{(i_5)}
+
{\bf 1}_{\{i_6=i_4\ne 0\}}
{\bf 1}_{\{j_6=j_4\}}
{\bf 1}_{\{i_3=i_5\ne 0\}}
{\bf 1}_{\{j_3=j_5\}}
\zeta_{j_1}^{(i_1)}
\zeta_{j_2}^{(i_2)}+
$$
$$
+
{\bf 1}_{\{i_6=i_4\ne 0\}}
{\bf 1}_{\{j_6=j_4\}}
{\bf 1}_{\{i_2=i_5\ne 0\}}
{\bf 1}_{\{j_2=j_5\}}
\zeta_{j_1}^{(i_1)}
\zeta_{j_3}^{(i_3)}
+
{\bf 1}_{\{i_6=i_4\ne 0\}}
{\bf 1}_{\{j_6=j_4\}}
{\bf 1}_{\{i_2=i_3\ne 0\}}
{\bf 1}_{\{j_2=j_3\}}
\zeta_{j_1}^{(i_1)}
\zeta_{j_5}^{(i_5)}+
$$
$$
+
{\bf 1}_{\{i_6=i_4\ne 0\}}
{\bf 1}_{\{j_6=j_4\}}
{\bf 1}_{\{i_1=i_5\ne 0\}}
{\bf 1}_{\{j_1=j_5\}}
\zeta_{j_2}^{(i_2)}
\zeta_{j_3}^{(i_3)}
+
{\bf 1}_{\{i_6=i_4\ne 0\}}
{\bf 1}_{\{j_6=j_4\}}
{\bf 1}_{\{i_1=i_3\ne 0\}}
{\bf 1}_{\{j_1=j_3\}}
\zeta_{j_2}^{(i_2)}
\zeta_{j_5}^{(i_5)}+
$$
$$
+
{\bf 1}_{\{i_6=i_4\ne 0\}}
{\bf 1}_{\{j_6=j_4\}}
{\bf 1}_{\{i_1=i_2\ne 0\}}
{\bf 1}_{\{j_1=j_2\}}
\zeta_{j_3}^{(i_3)}
\zeta_{j_5}^{(i_5)}
+
{\bf 1}_{\{i_6=i_5\ne 0\}}
{\bf 1}_{\{j_6=j_5\}}
{\bf 1}_{\{i_3=i_4\ne 0\}}
{\bf 1}_{\{j_3=j_4\}}
\zeta_{j_1}^{(i_1)}
\zeta_{j_2}^{(i_2)}+
$$
$$
+
{\bf 1}_{\{i_6=i_5\ne 0\}}
{\bf 1}_{\{j_6=j_5\}}
{\bf 1}_{\{i_2=i_4\ne 0\}}
{\bf 1}_{\{j_2=j_4\}}
\zeta_{j_1}^{(i_1)}
\zeta_{j_3}^{(i_3)}
+
{\bf 1}_{\{i_6=i_5\ne 0\}}
{\bf 1}_{\{j_6=j_5\}}
{\bf 1}_{\{i_2=i_3\ne 0\}}
{\bf 1}_{\{j_2=j_3\}}
\zeta_{j_1}^{(i_1)}
\zeta_{j_4}^{(i_4)}+
$$
$$
+
{\bf 1}_{\{i_6=i_5\ne 0\}}
{\bf 1}_{\{j_6=j_5\}}
{\bf 1}_{\{i_1=i_4\ne 0\}}
{\bf 1}_{\{j_1=j_4\}}
\zeta_{j_2}^{(i_2)}
\zeta_{j_3}^{(i_3)}
+
{\bf 1}_{\{i_6=i_5\ne 0\}}
{\bf 1}_{\{j_6=j_5\}}
{\bf 1}_{\{i_1=i_3\ne 0\}}
{\bf 1}_{\{j_1=j_3\}}
\zeta_{j_2}^{(i_2)}
\zeta_{j_4}^{(i_4)}+
$$
$$
+
{\bf 1}_{\{i_6=i_5\ne 0\}}
{\bf 1}_{\{j_6=j_5\}}
{\bf 1}_{\{i_1=i_2\ne 0\}}
{\bf 1}_{\{j_1=j_2\}}
\zeta_{j_3}^{(i_3)}
\zeta_{j_4}^{(i_4)}-
$$
$$
-
{\bf 1}_{\{i_6=i_1\ne 0\}}
{\bf 1}_{\{j_6=j_1\}}
{\bf 1}_{\{i_2=i_5\ne 0\}}
{\bf 1}_{\{j_2=j_5\}}
{\bf 1}_{\{i_3=i_4\ne 0\}}
{\bf 1}_{\{j_3=j_4\}}-
$$
$$
-
{\bf 1}_{\{i_6=i_1\ne 0\}}
{\bf 1}_{\{j_6=j_1\}}
{\bf 1}_{\{i_2=i_4\ne 0\}}
{\bf 1}_{\{j_2=j_4\}}
{\bf 1}_{\{i_3=i_5\ne 0\}}
{\bf 1}_{\{j_3=j_5\}}-
$$
$$
-
{\bf 1}_{\{i_6=i_1\ne 0\}}
{\bf 1}_{\{j_6=j_1\}}
{\bf 1}_{\{i_2=i_3\ne 0\}}
{\bf 1}_{\{j_2=j_3\}}
{\bf 1}_{\{i_4=i_5\ne 0\}}
{\bf 1}_{\{j_4=j_5\}}-
$$
$$
-
{\bf 1}_{\{i_6=i_2\ne 0\}}
{\bf 1}_{\{j_6=j_2\}}
{\bf 1}_{\{i_1=i_5\ne 0\}}
{\bf 1}_{\{j_1=j_5\}}
{\bf 1}_{\{i_3=i_4\ne 0\}}
{\bf 1}_{\{j_3=j_4\}}-
$$
$$
-
{\bf 1}_{\{i_6=i_2\ne 0\}}
{\bf 1}_{\{j_6=j_2\}}
{\bf 1}_{\{i_1=i_4\ne 0\}}
{\bf 1}_{\{j_1=j_4\}}
{\bf 1}_{\{i_3=i_5\ne 0\}}
{\bf 1}_{\{j_3=j_5\}}-
$$
$$
-
{\bf 1}_{\{i_6=i_2\ne 0\}}
{\bf 1}_{\{j_6=j_2\}}
{\bf 1}_{\{i_1=i_3\ne 0\}}
{\bf 1}_{\{j_1=j_3\}}
{\bf 1}_{\{i_4=i_5\ne 0\}}
{\bf 1}_{\{j_4=j_5\}}-
$$
$$
-
{\bf 1}_{\{i_6=i_3\ne 0\}}
{\bf 1}_{\{j_6=j_3\}}
{\bf 1}_{\{i_1=i_5\ne 0\}}
{\bf 1}_{\{j_1=j_5\}}
{\bf 1}_{\{i_2=i_4\ne 0\}}
{\bf 1}_{\{j_2=j_4\}}-
$$
$$
-
{\bf 1}_{\{i_6=i_3\ne 0\}}
{\bf 1}_{\{j_6=j_3\}}
{\bf 1}_{\{i_1=i_4\ne 0\}}
{\bf 1}_{\{j_1=j_4\}}
{\bf 1}_{\{i_2=i_5\ne 0\}}
{\bf 1}_{\{j_2=j_5\}}-
$$
$$
-
{\bf 1}_{\{i_3=i_6\ne 0\}}
{\bf 1}_{\{j_3=j_6\}}
{\bf 1}_{\{i_1=i_2\ne 0\}}
{\bf 1}_{\{j_1=j_2\}}
{\bf 1}_{\{i_4=i_5\ne 0\}}
{\bf 1}_{\{j_4=j_5\}}-
$$
$$
-
{\bf 1}_{\{i_6=i_4\ne 0\}}
{\bf 1}_{\{j_6=j_4\}}
{\bf 1}_{\{i_1=i_5\ne 0\}}
{\bf 1}_{\{j_1=j_5\}}
{\bf 1}_{\{i_2=i_3\ne 0\}}
{\bf 1}_{\{j_2=j_3\}}-
$$
$$
-
{\bf 1}_{\{i_6=i_4\ne 0\}}
{\bf 1}_{\{j_6=j_4\}}
{\bf 1}_{\{i_1=i_3\ne 0\}}
{\bf 1}_{\{j_1=j_3\}}
{\bf 1}_{\{i_2=i_5\ne 0\}}
{\bf 1}_{\{j_2=j_5\}}-
$$
$$
-
{\bf 1}_{\{i_6=i_4\ne 0\}}
{\bf 1}_{\{j_6=j_4\}}
{\bf 1}_{\{i_1=i_2\ne 0\}}
{\bf 1}_{\{j_1=j_2\}}
{\bf 1}_{\{i_3=i_5\ne 0\}}
{\bf 1}_{\{j_3=j_5\}}-
$$
$$
-
{\bf 1}_{\{i_6=i_5\ne 0\}}
{\bf 1}_{\{j_6=j_5\}}
{\bf 1}_{\{i_1=i_4\ne 0\}}
{\bf 1}_{\{j_1=j_4\}}
{\bf 1}_{\{i_2=i_3\ne 0\}}
{\bf 1}_{\{j_2=j_3\}}-
$$
$$
-
{\bf 1}_{\{i_6=i_5\ne 0\}}
{\bf 1}_{\{j_6=j_5\}}
{\bf 1}_{\{i_1=i_2\ne 0\}}
{\bf 1}_{\{j_1=j_2\}}
{\bf 1}_{\{i_3=i_4\ne 0\}}
{\bf 1}_{\{j_3=j_4\}}-
$$
\begin{equation}
\label{a6}
\Biggl.-
{\bf 1}_{\{i_6=i_5\ne 0\}}
{\bf 1}_{\{j_6=j_5\}}
{\bf 1}_{\{i_1=i_3\ne 0\}}
{\bf 1}_{\{j_1=j_3\}}
{\bf 1}_{\{i_2=i_4\ne 0\}}
{\bf 1}_{\{j_2=j_4\}}\Biggr),
\end{equation}

\vspace{6mm}
\noindent
where ${\bf 1}_A$ is the indicator of the set $A$.

For further consideration, let us 
consider the generalization of formulas (\ref{a1})--(\ref{a6})                 
for the case of an arbitrary multiplicity $k$ $(k\in\mathbb{N})$ of 
the iterated Ito stochastic integral $J[\psi^{(k)}]_{T,t}$ defined by (\ref{ito}).
In order to do this, let us
introduce some notations. 
Consider the unordered
set $\{1, 2, \ldots, k\}$ 
and separate it into two parts:
the first part consists of $r$ unordered 
pairs (sequence order of these pairs is also unimportant) and the 
second one consists of the 
remaining $k-2r$ numbers.
So, we have

\begin{equation}
\label{leto5007}
(\{
\underbrace{\{g_1, g_2\}, \ldots, 
\{g_{2r-1}, g_{2r}\}}_{\small{\hbox{part 1}}}
\},
\{\underbrace{q_1, \ldots, q_{k-2r}}_{\small{\hbox{part 2}}}
\}),
\end{equation}

\vspace{4mm}
\noindent
where 

\vspace{-2mm}
$$
\{g_1, g_2, \ldots, 
g_{2r-1}, g_{2r}, q_1, \ldots, q_{k-2r}\}=\{1, 2, \ldots, k\},
$$

\vspace{4mm}
\noindent
braces   
mean an unordered 
set, and pa\-ren\-the\-ses mean an ordered set.

We will say that (\ref{leto5007}) is a partition 
and consider the sum with respect to all possible
partitions

\begin{equation}
\label{leto5008}
\sum_{\stackrel{(\{\{g_1, g_2\}, \ldots, 
\{g_{2r-1}, g_{2r}\}\}, \{q_1, \ldots, q_{k-2r}\})}
{{}_{\{g_1, g_2, \ldots, 
g_{2r-1}, g_{2r}, q_1, \ldots, q_{k-2r}\}=\{1, 2, \ldots, k\}}}}
a_{g_1 g_2, \ldots, 
g_{2r-1} g_{2r}, q_1 \ldots q_{k-2r}}.
\end{equation}

\vspace{4mm}

Below there are several examples of sums in the form (\ref{leto5008})

\vspace{2mm}
$$
\sum_{\stackrel{(\{g_1, g_2\})}{{}_{\{g_1, g_2\}=\{1, 2\}}}}
a_{g_1 g_2}=a_{12},
$$

\vspace{3mm}
$$
\sum_{\stackrel{(\{\{g_1, g_2\}, \{g_3, g_4\}\})}
{{}_{\{g_1, g_2, g_3, g_4\}=\{1, 2, 3, 4\}}}}
a_{g_1 g_2, g_3 g_4}=a_{12,34} + a_{13,24} + a_{23,14},
$$

\vspace{3mm}
$$
\sum_{\stackrel{(\{g_1, g_2\}, \{q_1, q_{2}\})}
{{}_{\{g_1, g_2, q_1, q_{2}\}=\{1, 2, 3, 4\}}}}
a_{g_1 g_2, q_1 q_{2}}=
$$

$$
=a_{12,34}+a_{13,24}+a_{14,23}
+a_{23,14}+a_{24,13}+a_{34,12},
$$

\vspace{3mm}
$$
\sum_{\stackrel{(\{g_1, g_2\}, \{q_1, q_{2}, q_3\})}
{{}_{\{g_1, g_2, q_1, q_{2}, q_3\}=\{1, 2, 3, 4, 5\}}}}
a_{g_1 g_2, q_1 q_{2}q_3}
=
$$

$$
=a_{12,345}+a_{13,245}+a_{14,235}
+a_{15,234}+a_{23,145}+a_{24,135}+
$$
$$
+a_{25,134}+a_{34,125}+a_{35,124}+a_{45,123},
$$

\vspace{4mm}
$$
\sum_{\stackrel{(\{\{g_1, g_2\}, \{g_3, g_{4}\}\}, \{q_1\})}
{{}_{\{g_1, g_2, g_3, g_{4}, q_1\}=\{1, 2, 3, 4, 5\}}}}
a_{g_1 g_2, g_3 g_{4},q_1}
=
$$

$$
=
a_{12,34,5}+a_{13,24,5}+a_{14,23,5}+
a_{12,35,4}+a_{13,25,4}+a_{15,23,4}+
$$
$$
+a_{12,54,3}+a_{15,24,3}+a_{14,25,3}+a_{15,34,2}+a_{13,54,2}+a_{14,53,2}+
$$
$$
+
a_{52,34,1}+a_{53,24,1}+a_{54,23,1}.
$$

\vspace{5mm}

Now we can write (\ref{tyyy}) as

\vspace{1mm}

$$
J[\psi^{(k)}]_{T,t}=
\hbox{\vtop{\offinterlineskip\halign{
\hfil#\hfil\cr
{\rm l.i.m.}\cr
$\stackrel{}{{}_{p_1,\ldots,p_k\to \infty}}$\cr
}} }
\sum\limits_{j_1=0}^{p_1}\ldots
\sum\limits_{j_k=0}^{p_k}
C_{j_k\ldots j_1}\Biggl(
\prod_{l=1}^k\zeta_{j_l}^{(i_l)}+\sum\limits_{r=1}^{[k/2]}
(-1)^r \times
\Biggr.
$$

\vspace{3mm}
\begin{equation}
\label{leto6000hh}
\times
\sum_{\stackrel{(\{\{g_1, g_2\}, \ldots, 
\{g_{2r-1}, g_{2r}\}\}, \{q_1, \ldots, q_{k-2r}\})}
{{}_{\{g_1, g_2, \ldots, 
g_{2r-1}, g_{2r}, q_1, \ldots, q_{k-2r}\}=\{1, 2, \ldots, k\}}}}
\prod\limits_{s=1}^r
{\bf 1}_{\{i_{g_{{}_{2s-1}}}=~i_{g_{{}_{2s}}}\ne 0\}}
\Biggl.{\bf 1}_{\{j_{g_{{}_{2s-1}}}=~j_{g_{{}_{2s}}}\}}
\prod_{l=1}^{k-2r}\zeta_{j_{q_l}}^{(i_{q_l})}\Biggr),
\end{equation}

\vspace{5mm}
\noindent
where $[x]$ is an integer part of a real number $x;$
another notations are the same as in Theorem {\bf 1}.

\vspace{2mm}

In particular, from (\ref{leto6000hh}) for $k=5$ we obtain

\vspace{3mm}

$$
J[\psi^{(5)}]_{T,t}=
\hbox{\vtop{\offinterlineskip\halign{
\hfil#\hfil\cr
{\rm l.i.m.}\cr
$\stackrel{}{{}_{p_1,\ldots,p_5\to \infty}}$\cr
}} }\sum_{j_1=0}^{p_1}\ldots\sum_{j_5=0}^{p_5}
C_{j_5\ldots j_1}\Biggl(
\prod_{l=1}^5\zeta_{j_l}^{(i_l)}-\Biggr.
$$

\vspace{2mm}
$$
-
\sum\limits_{\stackrel{(\{g_1, g_2\}, \{q_1, q_{2}, q_3\})}
{{}_{\{g_1, g_2, q_{1}, q_{2}, q_3\}=\{1, 2, 3, 4, 5\}}}}
{\bf 1}_{\{i_{g_{{}_{1}}}=~i_{g_{{}_{2}}}\ne 0\}}
{\bf 1}_{\{j_{g_{{}_{1}}}=~j_{g_{{}_{2}}}\}}
\prod_{l=1}^{3}\zeta_{j_{q_l}}^{(i_{q_l})}+
$$

\vspace{2mm}
$$
+
\sum_{\stackrel{(\{\{g_1, g_2\}, 
\{g_{3}, g_{4}\}\}, \{q_1\})}
{{}_{\{g_1, g_2, g_{3}, g_{4}, q_1\}=\{1, 2, 3, 4, 5\}}}}
{\bf 1}_{\{i_{g_{{}_{1}}}=~i_{g_{{}_{2}}}\ne 0\}}
{\bf 1}_{\{j_{g_{{}_{1}}}=~j_{g_{{}_{2}}}\}}
\Biggl.{\bf 1}_{\{i_{g_{{}_{3}}}=~i_{g_{{}_{4}}}\ne 0\}}
{\bf 1}_{\{j_{g_{{}_{3}}}=~j_{g_{{}_{4}}}\}}
\zeta_{j_{q_1}}^{(i_{q_1})}\Biggr).
$$

\vspace{7mm}
\noindent
The last equality obviously agrees with
(\ref{a5}).

Let us consider a generalization of Theorem 1 for the case
of an arbitrary complete orthonormal systems  
of functions in the space $L_2([t,T])$ 
and $\psi_1(\tau),\ldots,\psi_k(\tau)\in L_2([t, T]).$

\vspace{2mm}

{\bf Theorem~2}\ \cite{2018a} (Sect.~1.11), \cite{arxiv-1} (Sect.~15).
{\it Suppose that
$\psi_1(\tau),\ldots,\psi_k(\tau)\in L_2([t, T])$ and
$\{\phi_j(x)\}_{j=0}^{\infty}$ is an arbitrary complete orthonormal system  
of functions in the space $L_2([t,T]).$
Then the following expansion

\vspace{1mm}
$$
J[\psi^{(k)}]_{T,t}=
\hbox{\vtop{\offinterlineskip\halign{
\hfil#\hfil\cr
{\rm l.i.m.}\cr
$\stackrel{}{{}_{p_1,\ldots,p_k\to \infty}}$\cr
}} }
\sum\limits_{j_1=0}^{p_1}\ldots
\sum\limits_{j_k=0}^{p_k}
C_{j_k\ldots j_1}\Biggl(
\prod_{l=1}^k\zeta_{j_l}^{(i_l)}+\sum\limits_{r=1}^{[k/2]}
(-1)^r \times
\Biggr.
$$

\vspace{2mm}
\begin{equation}
\label{leto6000}
\times
\sum_{\stackrel{(\{\{g_1, g_2\}, \ldots, 
\{g_{2r-1}, g_{2r}\}\}, \{q_1, \ldots, q_{k-2r}\})}
{{}_{\{g_1, g_2, \ldots, 
g_{2r-1}, g_{2r}, q_1, \ldots, q_{k-2r}\}=\{1, 2, \ldots, k\}}}}
\prod\limits_{s=1}^r
{\bf 1}_{\{i_{g_{{}_{2s-1}}}=~i_{g_{{}_{2s}}}\ne 0\}}
\Biggl.{\bf 1}_{\{j_{g_{{}_{2s-1}}}=~j_{g_{{}_{2s}}}\}}
\prod_{l=1}^{k-2r}\zeta_{j_{q_l}}^{(i_{q_l})}\Biggr)
\end{equation}

\vspace{6mm}
\noindent
con\-verg\-ing in the mean-square sense is valid,
where $[x]$ is an integer part of a real number $x;$
another notations are the same as in Theorem~{\rm 1}.}

\vspace{2mm}

It should be noted that an analogue of Theorem 2 was considered 
in \cite{Rybakov1000}. 
Note that we use another notations 
\cite{2018a} (Sect.~1.11), \cite{arxiv-1} (Sect.~15)
in comparison with \cite{Rybakov1000}.
Moreover, the proof of an analogue of Theorem 2
from \cite{Rybakov1000} is somewhat different from the proof given in 
\cite{2018a} (Sect.~1.11), \cite{arxiv-1} (Sect.~15).

\vspace{5mm}

\section{Expansions of Iterated Stratonovich Stochastic Integrals of Multiplicities 2 to 6}

\vspace{5mm}

In a number of works of the author 
\cite{2010-2}-\cite{2018axx}, \cite{arxiv-2}
Theorems 1, 2 have been adapted for the iterated  Stratonovich stochastic integrals
(\ref{str}) of multiplicities 2 to 6.
Let us first present some old results as the following theorem.

\vspace{2mm}

{\bf Theorem 3} \cite{2010-2}-\cite{2018axx}, \cite{arxiv-2}.\
{\it Suppose that 
$\{\phi_j(x)\}_{j=0}^{\infty}$ is a complete orthonormal system of 
Legendre polynomials or trigonometric functions in the space $L_2([t, T]).$
At the same time $\psi_2(\tau)$ is a continuously differentiable 
function on $[t, T]$ and $\psi_1(\tau), \psi_3(\tau)$ are twice 
continuously differentiable functions on $[t, T]$. Then

\begin{equation}
\label{a}
J^{*}[\psi^{(2)}]_{T,t}=
\hbox{\vtop{\offinterlineskip\halign{
\hfil#\hfil\cr
{\rm l.i.m.}\cr
$\stackrel{}{{}_{p_1,p_2\to \infty}}$\cr
}} }\sum_{j_1=0}^{p_1}\sum_{j_2=0}^{p_2}
C_{j_2j_1}\zeta_{j_1}^{(i_1)}\zeta_{j_2}^{(i_2)}\ \ \ (i_1,i_2=1,\ldots,m),
\end{equation}

\vspace{1mm}
\begin{equation}
\label{feto19000ab}
J^{*}[\psi^{(3)}]_{T,t}=
\hbox{\vtop{\offinterlineskip\halign{
\hfil#\hfil\cr
{\rm l.i.m.}\cr
$\stackrel{}{{}_{p_1,p_2,p_3\to \infty}}$\cr
}} }\sum_{j_1=0}^{p_1}\sum_{j_2=0}^{p_2}\sum_{j_3=0}^{p_3}
C_{j_3 j_2 j_1}\zeta_{j_1}^{(i_1)}\zeta_{j_2}^{(i_2)}\zeta_{j_3}^{(i_3)}\ \ \
(i_1,i_2,i_3=0, 1,\ldots,m),
\end{equation}

\vspace{1mm}
\begin{equation}
\label{feto19000a}
J^{*}[\psi^{(3)}]_{T,t}=
\hbox{\vtop{\offinterlineskip\halign{
\hfil#\hfil\cr
{\rm l.i.m.}\cr
$\stackrel{}{{}_{p\to \infty}}$\cr
}} }
\sum\limits_{j_1,j_2,j_3=0}^{p}
C_{j_3 j_2 j_1}\zeta_{j_1}^{(i_1)}\zeta_{j_2}^{(i_2)}\zeta_{j_3}^{(i_3)}\ \ \
(i_1,i_2,i_3=1,\ldots,m),
\end{equation}

\vspace{1mm}
\begin{equation}
\label{uu}
J^{*}[\psi^{(4)}]_{T,t}=
\hbox{\vtop{\offinterlineskip\halign{
\hfil#\hfil\cr
{\rm l.i.m.}\cr
$\stackrel{}{{}_{p\to \infty}}$\cr
}} }
\sum\limits_{j_1,j_2,j_3,j_4=0}^{p}
C_{j_4 j_3 j_2 j_1}\zeta_{j_1}^{(i_1)}
\zeta_{j_2}^{(i_2)}\zeta_{j_3}^{(i_3)}\zeta_{j_4}^{(i_4)}\ \ \
(i_1,i_2,i_3,i_4=0, 1,\ldots,m),
\end{equation}

\vspace{5mm}
\noindent
where $J^{*}[\psi^{(k)}]_{T,t}$ is defined by {\rm (\ref{str})} and
$\psi_l(\tau)\equiv 1$ $(l=1,\ldots,4)$ in {\rm (\ref{feto19000ab})}, 
{\rm (\ref{uu});} another notations are the same as in Theorems {\rm 1, 2.}
}

\vspace{2mm}

Recently, a new approach to the expansion and mean-square 
approximation of iterated Stratonovich stochastic integrals has been obtained
\cite{2018a} (Sect.~2.10--2.16), \cite{arxiv-2} (Sect.~13--19), 
\cite{arxiv-4} (Sect.~7--13), \cite{arxiv-5} (Sect.~5--11),
\cite{new-art-1-xxy} (Sect.~4--9), \cite{new-art-1xxys}.
Let us formulate four theorems that were obtained using this approach.

\vspace{2mm}

{\bf Theorem 4}\ \cite{2018a}, \cite{arxiv-2}, \cite{arxiv-4}, \cite{arxiv-5}, \cite{new-art-1-xxy}.\
{\it Suppose 
that $\{\phi_j(x)\}_{j=0}^{\infty}$ is a complete orthonormal system of 
Legendre polynomials or trigonometric functions in the space $L_2([t, T]).$
Furthermore, let $\psi_1(\tau), \psi_2(\tau),$ $\psi_3(\tau)$ are continuously dif\-ferentiable 
nonrandom functions on $[t, T].$ 
Then, for the 
iterated Stra\-to\-no\-vich stochastic integral of third multiplicity

$$
J^{*}[\psi^{(3)}]_{T,t}={\int\limits_t^{*}}^T\psi_3(t_3)
{\int\limits_t^{*}}^{t_3}\psi_2(t_2)
{\int\limits_t^{*}}^{t_2}\psi_1(t_1)
d{\bf w}_{t_1}^{(i_1)}
d{\bf w}_{t_2}^{(i_2)}d{\bf w}_{t_3}^{(i_3)}\ \ \ (i_1,i_2,i_3=0,1,\ldots,m)
$$

\vspace{4mm}
\noindent
the following 
relations

\vspace{-1mm}
\begin{equation}
\label{fin1}
J^{*}[\psi^{(3)}]_{T,t}
=\hbox{\vtop{\offinterlineskip\halign{
\hfil#\hfil\cr
{\rm l.i.m.}\cr
$\stackrel{}{{}_{p\to \infty}}$\cr
}} }
\sum\limits_{j_1, j_2, j_3=0}^{p}
C_{j_3 j_2 j_1}\zeta_{j_1}^{(i_1)}\zeta_{j_2}^{(i_2)}\zeta_{j_3}^{(i_3)},
\end{equation}

\vspace{3mm}
\begin{equation}
\label{fin2}
{\sf M}\left\{\left(
J^{*}[\psi^{(3)}]_{T,t}-
\sum\limits_{j_1, j_2, j_3=0}^{p}
C_{j_3 j_2 j_1}\zeta_{j_1}^{(i_1)}\zeta_{j_2}^{(i_2)}\zeta_{j_3}^{(i_3)}\right)^2\right\}
\le \frac{C}{p}
\end{equation}

\vspace{5mm}
\noindent
are fulfilled, where $i_1, i_2, i_3=0,1,\ldots,m$ in {\rm (\ref{fin1})} and 
$i_1, i_2, i_3=1,\ldots,m$ in {\rm (\ref{fin2})},
constant $C$ is independent of $p,$

$$
C_{j_3 j_2 j_1}=\int\limits_t^T\psi_3(t_3)\phi_{j_3}(t_3)
\int\limits_t^{t_3}\psi_2(t_2)\phi_{j_2}(t_2)
\int\limits_t^{t_2}\psi_1(t_1)\phi_{j_1}(t_1)dt_1dt_2dt_3
$$

\vspace{4mm}
\noindent
and
$$
\zeta_{j}^{(i)}=
\int\limits_t^T \phi_{j}(\tau) d{\bf f}_{\tau}^{(i)}
$$ 

\vspace{2mm}
\noindent
are independent standard Gaussian random variables for various 
$i$ or $j$ {\rm (}in the case when $i\ne 0${\rm );} 
another notations are the same as in Theorems~{\rm 1, 2}.}

\vspace{2mm}

{\bf Theorem 5}\ \cite{2018a}, \cite{arxiv-2}, \cite{arxiv-4}, \cite{arxiv-5}, \cite{new-art-1-xxy}.\ 
{\it Let
$\{\phi_j(x)\}_{j=0}^{\infty}$ be a complete orthonormal system of 
Legendre polynomials or trigonometric functions in the space $L_2([t, T]).$
Furthermore, let $\psi_1(\tau), \ldots, \psi_4(\tau)$ be continuously dif\-ferentiable 
nonrandom functions on $[t, T].$ 
Then, for the 
iterated Stra\-to\-no\-vich stochastic integral of fourth multiplicity

\begin{equation}
\label{fin0}
J^{*}[\psi^{(4)}]_{T,t}={\int\limits_t^{*}}^T\psi_4(t_4)
{\int\limits_t^{*}}^{t_4}\psi_3(t_3)
{\int\limits_t^{*}}^{t_3}\psi_2(t_2)
{\int\limits_t^{*}}^{t_2}\psi_1(t_1)
d{\bf w}_{t_1}^{(i_1)}
d{\bf w}_{t_2}^{(i_2)}d{\bf w}_{t_3}^{(i_3)}d{\bf w}_{t_4}^{(i_4)}
\end{equation}

\vspace{4mm}
\noindent
the following 
relations

\begin{equation}
\label{fin3}
J^{*}[\psi^{(4)}]_{T,t}
=\hbox{\vtop{\offinterlineskip\halign{
\hfil#\hfil\cr
{\rm l.i.m.}\cr
$\stackrel{}{{}_{p\to \infty}}$\cr
}} }
\sum\limits_{j_1, j_2, j_3,j_4=0}^{p}
C_{j_4j_3 j_2 j_1}\zeta_{j_1}^{(i_1)}\zeta_{j_2}^{(i_2)}\zeta_{j_3}^{(i_3)}\zeta_{j_4}^{(i_4)},
\end{equation}

\vspace{3mm}

\begin{equation}
\label{fin4}
{\sf M}\left\{\left(
J^{*}[\psi^{(4)}]_{T,t}-
\sum\limits_{j_1, j_2, j_3, j_4=0}^{p}
C_{j_4 j_3 j_2 j_1}\zeta_{j_1}^{(i_1)}\zeta_{j_2}^{(i_2)}\zeta_{j_3}^{(i_3)}
\zeta_{j_4}^{(i_4)}
\right)^2\right\}
\le \frac{C}{p^{1-\varepsilon}}
\end{equation}

\vspace{5mm}
\noindent
are fulfilled, where $i_1, \ldots , i_4=0,1,\ldots,m$ in {\rm (\ref{fin0}),} {\rm (\ref{fin3})} 
and $i_1, \ldots, i_4=1,\ldots,m$ in {\rm (\ref{fin4}),}
constant $C$ does not depend on $p,$
$\varepsilon$ is an arbitrary
small positive real number 
for the case of complete orthonormal system of 
Legendre polynomials in the space $L_2([t, T])$
and $\varepsilon=0$ for the case of
complete orthonormal system of 
trigonometric functions in the space $L_2([t, T]),$

$$
C_{j_4 j_3 j_2 j_1}=
$$

$$
=
\int\limits_t^T\psi_4(t_4)\phi_{j_4}(t_4)
\int\limits_t^{t_4}\psi_3(t_3)\phi_{j_3}(t_3)
\int\limits_t^{t_3}\psi_2(t_2)\phi_{j_2}(t_2)
\int\limits_t^{t_2}\psi_1(t_1)\phi_{j_1}(t_1)dt_1dt_2dt_3dt_4;
$$

\vspace{4mm}
\noindent
another notations are the same as in Theorem~{\rm 4}.}

\vspace{2mm}

{\bf Theorem 6}\ \cite{2018a}, \cite{arxiv-2}, \cite{arxiv-4}, \cite{arxiv-5}, \cite{new-art-1-xxy}.\
{\it Assume 
that $\{\phi_j(x)\}_{j=0}^{\infty}$ is a complete orthonormal system of 
Legendre polynomials or trigonometric functions in the space $L_2([t, T])$
and $\psi_1(\tau), \ldots, \psi_5(\tau)$ are continuously dif\-ferentiable 
nonrandom functions on $[t, T].$ 
Then, for the 
iterated Stra\-to\-no\-vich stochastic integral of fifth multiplicity

\begin{equation}
\label{fin7}
J^{*}[\psi^{(5)}]_{T,t}={\int\limits_t^{*}}^T\psi_5(t_5)
\ldots
{\int\limits_t^{*}}^{t_2}\psi_1(t_1)
d{\bf w}_{t_1}^{(i_1)}
\ldots d{\bf w}_{t_5}^{(i_5)}
\end{equation}

\vspace{4mm}
\noindent
the following 
relations

\begin{equation}
\label{fin8}
J^{*}[\psi^{(5)}]_{T,t}
=\hbox{\vtop{\offinterlineskip\halign{
\hfil#\hfil\cr
{\rm l.i.m.}\cr
$\stackrel{}{{}_{p\to \infty}}$\cr
}} }
\sum\limits_{j_1,\ldots,j_5=0}^{p}
C_{j_5 \ldots j_1}\zeta_{j_1}^{(i_1)}\ldots \zeta_{j_5}^{(i_5)},
\end{equation}

\vspace{3mm}

\begin{equation}
\label{fin9}
{\sf M}\left\{\left(
J^{*}[\psi^{(5)}]_{T,t}-
\sum\limits_{j_1, \ldots, j_5=0}^{p}
C_{j_5 \ldots j_1}\zeta_{j_1}^{(i_1)}\ldots
\zeta_{j_5}^{(i_5)}
\right)^2\right\}
\le \frac{C}{p^{1-\varepsilon}}
\end{equation}

\vspace{5mm}
\noindent
are fulfilled, where $i_1, \ldots , i_5=0,1,\ldots,m$ in {\rm (\ref{fin7}),} {\rm (\ref{fin8})} 
and $i_1, \ldots, i_5=1,\ldots,m$ in {\rm (\ref{fin9}),}
constant $C$ is independent of $p,$
$\varepsilon$ is an arbitrary
small positive real number 
for the case of complete orthonormal system of 
Legendre polynomials in the space $L_2([t, T])$
and $\varepsilon=0$ for the case of
complete orthonormal system of 
trigonometric functions in the space $L_2([t, T]),$

$$
C_{j_5 \ldots j_1}=
\int\limits_t^T\psi_5(t_5)\phi_{j_5}(t_5)\ldots
\int\limits_t^{t_2}\psi_1(t_1)\phi_{j_1}(t_1)dt_1\ldots dt_5;
$$

\vspace{3mm}
\noindent
another notations are the same as in Theorems~{\rm 4, 5}.}

\vspace{2mm}

{\bf Theorem 7}\ \cite{2018a}, \cite{arxiv-2}, \cite{arxiv-4}, \cite{arxiv-5}, \cite{new-art-1xxys}.\
{\it Suppose that 
$\{\phi_j(x)\}_{j=0}^{\infty}$ is a complete orthonormal system of 
Legendre polynomials or trigonometric functions in the space $L_2([t, T]).$
Then, for the 
iterated Stratonovich stochastic integral of sixth multiplicity

\begin{equation}
\label{after10001qu1}
J_{T,t}^{*(i_1\ldots i_6)}={\int\limits_t^{*}}^T
\ldots
{\int\limits_t^{*}}^{t_2}
d{\bf w}_{t_1}^{(i_1)}
\ldots d{\bf w}_{t_6}^{(i_6)}
\end{equation}

\vspace{3mm}
\noindent
the following 
expansion 

\vspace{-1mm}
$$
J_{T,t}^{*(i_1\ldots i_6)}
=\hbox{\vtop{\offinterlineskip\halign{
\hfil#\hfil\cr
{\rm l.i.m.}\cr
$\stackrel{}{{}_{p\to \infty}}$\cr
}} }
\sum\limits_{j_1, \ldots, j_6=0}^{p}
C_{j_6 \ldots j_1}\zeta_{j_1}^{(i_1)}\ldots
\zeta_{j_6}^{(i_6)}
$$

\vspace{4mm}
\noindent
that converges in the mean-square sense is valid, where
$i_1, \ldots, i_6=0, 1,\ldots,m,$

$$
C_{j_6 \ldots j_1}=
\int\limits_t^T\phi_{j_6}(t_6)\ldots
\int\limits_t^{t_2}\phi_{j_1}(t_1)dt_1\ldots dt_6;
$$

\vspace{3mm}
\noindent
another notations are the same as in Theorems~{\rm 4--6}.}

\vspace{2mm}

The results of Theorems~3--7 were 
developed in \cite{2018a} (Chapter~2), \cite{arxiv-2}, \cite{arxiv-4}, \cite{arxiv-5}.
In particular, analogues of Theorem~7 for iterated Stratonovich stochastic
integrals of multiplicities 7 and 8 were obtained in \cite{2018a} (Sect.~2.36, 2.37).
In addition, the variants of Thorems 3--7 
were obtained
for the case when $\{\phi_j(x)\}_{j=0}^{\infty}$ is an arbitrary complete orthonormal system
of functions in $L_2([t, T])$ \cite{2018a} (Sect.~2.1.4, 2.23, 2.24, 2.31--2.34),
\cite{arxiv-2}, \cite{arxiv-4}, \cite{arxiv-5}.

\vspace{5mm}

\section{Exact Calculation of the Mean-Square Error in Theorems~1, 2}

\vspace{5mm}

Theorems 1 and 2 allow us to accurately calculate 
the mean-square 
approximation error for iterated Ito stochastic integrals
(see Theorem 8 below).

Assume that $J[\psi^{(k)}]_{T,t}^{p_1 \ldots p_k}$ is the approximation 
of (\ref{ito}), which is
the expression on the right-hand side of (\ref{leto6000}) before passing to the limit

\vspace{1mm}
$$
J[\psi^{(k)}]_{T,t}^{p_1 \ldots p_k}=
\sum\limits_{j_1=0}^{p_1}\ldots
\sum\limits_{j_k=0}^{p_k}
C_{j_k\ldots j_1}\Biggl(
\prod_{l=1}^k\zeta_{j_l}^{(i_l)}+\sum\limits_{r=1}^{[k/2]}
(-1)^r \times
\Biggr.
$$

\vspace{4mm}
$$
\times
\sum_{\stackrel{(\{\{g_1, g_2\}, \ldots, 
\{g_{2r-1}, g_{2r}\}\}, \{q_1, \ldots, q_{k-2r}\})}
{{}_{\{g_1, g_2, \ldots, 
g_{2r-1}, g_{2r}, q_1, \ldots, q_{k-2r}\}=\{1, 2, \ldots, k\}}}}
\prod\limits_{s=1}^r
{\bf 1}_{\{i_{g_{{}_{2s-1}}}=~i_{g_{{}_{2s}}}\ne 0\}}
\Biggl.{\bf 1}_{\{j_{g_{{}_{2s-1}}}=~j_{g_{{}_{2s}}}\}}
\prod_{l=1}^{k-2r}\zeta_{j_{q_l}}^{(i_{q_l})}\Biggr),
$$

\vspace{6mm}
\noindent
where $[x]$ is an integer part of a real number $x;$
another notations are the same as in Theorems~{\rm 1, 2}.

Let us denote

$$
E_k^{p_1,\ldots,p_k}\stackrel{{\rm def}}
{=}{\sf M}\left\{\left(J[\psi^{(k)}]_{T,t}-
J[\psi^{(k)}]_{T,t}^{p_1,\ldots,p_k}\right)^2\right\},
$$

\vspace{3mm}
$$
E_k^{p_1,\ldots,p_k}\stackrel{{\rm def}}{=}E_k^p\stackrel\ \ \hbox{if}\ \ 
p_1=\ldots=p_k=p,
$$

\vspace{2mm}
$$
I_k\stackrel{{\rm def}}{=}\left\Vert K\right\Vert^2_{L_2([t,T]^k)}=\int\limits_{[t,T]^k}
K^2(t_1,\ldots,t_k)dt_1\ldots dt_k.
$$

\vspace{4mm}

In \cite{2006}-\cite{2018axx}, \cite{arxiv-1} it was shown that

\begin{equation}
\label{star00011}
E_k^{p_1,\ldots,p_k}\le k!\left(I_k-\sum_{j_1=0}^{p_1}\ldots
\sum_{j_k=0}^{p_k}C^2_{j_k\ldots j_1}\right)
\end{equation}

\vspace{4mm}
\noindent
if $i_1,\ldots,i_k=1,\ldots,m$ and $0<T-t<\infty$ or 
$i_1,\ldots,i_k=0, 1,\ldots,m$ and $0<T-t<1.$

Moreover,    
in \cite{2018a} (Sect.~1.1.9, 1.11, 1.12), \cite{arxiv-1} (Sect.~6, 15, 16)
the following estimate 

\vspace{1mm}
\begin{equation}
\label{99999}
{\sf M}\left\{\left(J[\psi^{(k)}]_{T,t}-
J[\psi^{(k)}]_{T,t}^{p_1,\ldots,p_k}\right)^{2n}\right\}\le
(k!)^{n}(2n-1)^{nk}\
\left(I_k-\sum_{j_1=0}^{p_1}\ldots
\sum_{j_k=0}^{p_k}C^2_{j_k\ldots j_1}\right)^n
\end{equation}

\vspace{5mm}
\noindent
is obtained, where $n\in \mathbb{N}$.

The value $E_k^{p}$
can be calculated exactly.

\vspace{2mm}

{\bf Theorem 8} \cite{2018a} (Sect.~1.12), \cite{arxiv-6} (Sect.~6).
{\it Suppose that $\{\phi_j(x)\}_{j=0}^{\infty}$ 
is an arbitrary complete orthonormal system  
of functions in the space $L_2([t,T])$ and
$\psi_1(\tau),\ldots,\psi_k(\tau)\in L_2([t, T]),$  $i_1,\ldots, i_k=1,\ldots,m$.
Then

\begin{equation}
\label{tttr11}
E_k^p=I_k- \sum_{j_1,\ldots, j_k=0}^{p}
C_{j_k\ldots j_1}
{\sf M}\left\{J[\psi^{(k)}]_{T,t}
\sum\limits_{(j_1,\ldots,j_k)}
\int\limits_t^T \phi_{j_k}(t_k)
\ldots
\int\limits_t^{t_{2}}\phi_{j_{1}}(t_{1})
d{\bf f}_{t_1}^{(i_1)}\ldots
d{\bf f}_{t_k}^{(i_k)}\right\},
\end{equation}

\vspace{5mm}
\noindent
where $i_1,\ldots,i_k = 1,\ldots,m;$
the expression 

\vspace{-1mm}
$$
\sum\limits_{(j_1,\ldots,j_k)}
$$ 

\vspace{3mm}
\noindent
means the sum with respect to all
possible permutations 
$(j_1,\ldots,j_k)$. At the same time if 
$j_r$ swapped with $j_q$ in the permutation $(j_1,\ldots,j_k),$
then $i_r$ swapped with $i_q$ in the permutation
$(i_1,\ldots,i_k);$
another notations are the same as in Theorems {\rm 1, 2.}
}

Note that 

$$
{\sf M}\left\{J[\psi^{(k)}]_{T,t}
\int\limits_t^T \phi_{j_k}(t_k)
\ldots
\int\limits_t^{t_{2}}\phi_{j_{1}}(t_{1})
d{\bf f}_{t_1}^{(i_1)}\ldots
d{\bf f}_{t_k}^{(i_k)}\right\}=C_{j_k\ldots j_1}.
$$

\vspace{5mm}

Then from Theorem 8 for pairwise different $i_1,\ldots,i_k$ 
and for $i_1=\ldots=i_k$
we obtain

$$
E_k^p= I_k- \sum_{j_1,\ldots,j_k=0}^{p}
C_{j_k\ldots j_1}^2,
$$

\vspace{2mm}
$$ 
E_k^p= I_k - \sum_{j_1,\ldots,j_k=0}^{p}
C_{j_k\ldots j_1}\Biggl(\sum\limits_{(j_1,\ldots,j_k)}
C_{j_k\ldots j_1}\Biggr).
$$

\vspace{6mm}

Consider some examples of the application of Theorem 8
$(i_1, i_2 ,i_3=1,\ldots,m)$

\vspace{1mm}
$$
E_2^p
=I_2
-\sum_{j_1,j_2=0}^p
C_{j_2j_1}^2-
\sum_{j_1,j_2=0}^p
C_{j_2j_1}C_{j_1j_2}\ \ \ (i_1=i_2),
$$

\vspace{3mm}
$$
E_3^p=I_3
-\sum_{j_3,j_2,j_1=0}^p C_{j_3j_2j_1}^2-
\sum_{j_3,j_2,j_1=0}^p C_{j_3j_1j_2}C_{j_3j_2j_1}\ \ \ (i_1=i_2\ne i_3),
$$

\vspace{3mm}
$$
E_3^p=I_3-
\sum_{j_3,j_2,j_1=0}^p C_{j_3j_2j_1}^2-
\sum_{j_3,j_2,j_1=0}^p C_{j_2j_3j_1}C_{j_3j_2j_1}\ \ \ (i_1\ne i_2=i_3),
$$

\vspace{3mm}
$$
E_3^p=I_3
-\sum_{j_3,j_2,j_1=0}^p C_{j_3j_2j_1}^2-
\sum_{j_3,j_2,j_1=0}^p C_{j_3j_2j_1}C_{j_1j_2j_3}\ \ \ (i_1=i_3\ne i_2),
$$

\vspace{3mm}
$$
E^p_4 = I_4 - \sum_{j_1,\ldots,j_4=0}^{p}
C_{j_4\ldots j_1}\Biggl(\sum\limits_{(j_1,j_2)}
C_{j_4\ldots j_1}\Biggr)\ \ \ (i_1=i_2\ne i_3, i_4;\ i_3\ne i_4),
$$

\vspace{3mm}
$$
E^p_4 = I_4 - \sum_{j_1,\ldots,j_4=0}^{p}
C_{j_4\ldots j_1}\Biggl(\sum\limits_{(j_1,j_3)}
C_{j_4\ldots j_1}\Biggr)\ \ \ (i_1=i_3\ne i_2, i_4;\ i_2\ne i_4),
$$

\vspace{3mm}
$$
E^p_4 = I_4 - \sum_{j_1,\ldots,j_4=0}^{p}
C_{j_4\ldots j_1}\Biggl(\sum\limits_{(j_1,j_4)}
C_{j_4\ldots j_1}\Biggr)\ \ \ (i_1=i_4\ne i_2, i_3;\ i_2\ne i_3),
$$

\vspace{3mm}
$$
E^p_4 = I_4 - \sum_{j_1,\ldots,j_4=0}^{p}
C_{j_4\ldots j_1}\Biggl(\sum\limits_{(j_2,j_3)}
C_{j_4\ldots j_1}\Biggr)\ \ \ (i_2=i_3\ne i_1, i_4;\ i_1\ne i_4),
$$

\vspace{3mm}
$$
E^p_4 = I_4 - \sum_{j_1,\ldots,j_4=0}^{p}
C_{j_4\ldots j_1}\Biggl(\sum\limits_{(j_2,j_4)}
C_{j_4\ldots j_1}\Biggr)\ \ \ (i_2=i_4\ne i_1, i_3;\ i_1\ne i_3),
$$

\vspace{3mm}
$$
E^p_4 = I_4 - \sum_{j_1,\ldots,j_4=0}^{p}
C_{j_4\ldots j_1}\Biggl(\sum\limits_{(j_3,j_4)}
C_{j_4\ldots j_1}\Biggr)\ \ \ (i_3=i_4\ne i_1, i_2;\ i_1\ne i_2),
$$

\vspace{3mm}
$$
E_4^p = I_4 -
\sum_{j_1,\ldots,j_4=0}^{p}
C_{j_4\ldots j_1}\Biggl(\sum\limits_{(j_1,j_2,j_3)}
C_{j_4\ldots j_1}\Biggr)\ \ \ (i_1=i_2=i_3\ne i_4),
$$

\vspace{3mm}
$$
E_4^p = I_4 -
 \sum_{j_1,\ldots,j_4=0}^{p}
C_{j_4\ldots j_1}\Biggl(\sum\limits_{(j_2,j_3,j_4)}
C_{j_4\ldots j_1}\Biggr)\ \ \ (i_2=i_3=i_4\ne i_1),
$$

\vspace{3mm}
$$
E_4^p = I_4 -
 \sum_{j_1,\ldots,j_4=0}^{p}
C_{j_4\ldots j_1}\Biggl(\sum\limits_{(j_1,j_2,j_4)}
C_{j_4\ldots j_1}\Biggr)\ \ \ (i_1=i_2=i_4\ne i_3),
$$

\vspace{3mm}
$$
E_4^p = I_4 -
 \sum_{j_1,\ldots,j_4=0}^{p}
C_{j_4\ldots j_1}\Biggl(\sum\limits_{(j_1,j_3,j_4)}
C_{j_4\ldots j_1}\Biggr)\ \ \ (i_1=i_3=i_4\ne i_2),
$$

\vspace{3mm}
$$
E^p_4 = I_4 - \sum_{j_1,\ldots,j_4=0}^{p}
C_{j_4\ldots j_1}\Biggl(\sum\limits_{(j_1,j_2)}\Biggl(
\sum\limits_{(j_3,j_4)}
C_{j_4\ldots j_1}\Biggr)\Biggr)\ \ \ (i_1=i_2\ne i_3=i_4),
$$

\vspace{3mm}
$$
E^p_4 = I_4 - \sum_{j_1,\ldots,j_4=0}^{p}
C_{j_4\ldots j_1}\Biggl(\sum\limits_{(j_1,j_3)}\Biggl(
\sum\limits_{(j_2,j_4)}
C_{j_4\ldots j_1}\Biggr)\Biggr)\ \ \ (i_1=i_3\ne i_2=i_4),
$$

\vspace{3mm}
$$
E^p_4 = I_4 - \sum_{j_1,\ldots,j_4=0}^{p}
C_{j_4\ldots j_1}\Biggl(\sum\limits_{(j_1,j_4)}\Biggl(
\sum\limits_{(j_2,j_3)}
C_{j_4\ldots j_1}\Biggr)\Biggr)\ \ \ (i_1=i_4\ne i_2=i_3),
$$

\vspace{4mm}
$$
E_5^p = I_5 - \sum_{j_1,\ldots,j_5=0}^{p}
C_{j_5\ldots j_1}\Biggl(\sum\limits_{(j_2,j_4)}\Biggl(
\sum\limits_{(j_3,j_5)}
C_{j_5\ldots j_1}\Biggr)\Biggr)\ \ \ (i_1\ne i_2=i_4\ne i_3=i_5\ne i_1).
$$

\vspace{5mm}

\section{Some Technical Problems of the Milstein Approach}

\vspace{5mm}

Let us denote

\vspace{-2mm}
\begin{equation}
\label{k1001}
I_{(l_1\ldots l_k)T,t}^{*(i_1\ldots i_k)}
={\int\limits_t^{*}}^T(t-t_k)^{l_k} \ldots {\int\limits_t^{*}}^{t_2}
(t-t_1)^{l_1} d{\bf f}_{t_1}^{(i_1)}\ldots
d{\bf f}_{t_k}^{(i_k)},
\end{equation}

\vspace{4mm}
\noindent
where $i_1,\ldots, i_k=1,\dots,m;$ $l_1,\ldots,l_k=0, 1,\ldots.$

Consider the Milstein expansions for the simplest 
iterated Stratonovich stochastic integrals (\ref{k1001})

\begin{equation}
\label{41}
I_{(0)T,t}^{*(i_1)}=\sqrt{T-t}\zeta_0^{(i_1)},
\end{equation}

\vspace{2mm}
\begin{equation}
\label{42}
I_{(1)T,t}^{*(i_1)}=-\frac{{(T-t)}^{3/2}}{2}
\Biggl(\zeta_0^{(i_1)}-\frac{\sqrt{2}}{\pi}\sum_{r=1}^{\infty}
\frac{1}{r}\zeta_{2r-1}^{(i_1)}\Biggr),
\end{equation}

\vspace{4mm}
$$
I_{(00)T,t}^{*(i_1 i_2)}=\frac{1}{2}(T-t)\Biggl(
\zeta_{0}^{(i_1)}\zeta_{0}^{(i_2)}
+\frac{1}{\pi}
\sum_{r=1}^{\infty}\frac{1}{r}\left(
\zeta_{2r}^{(i_1)}\zeta_{2r-1}^{(i_2)}-
\zeta_{2r-1}^{(i_1)}\zeta_{2r}^{(i_2)}+
\right.\Biggr.
$$

\vspace{1mm}
\begin{equation}
\label{43}
+\left.\sqrt{2}\left(\zeta_{2r-1}^{(i_1)}\zeta_{0}^{(i_2)}-
\zeta_{0}^{(i_1)}\zeta_{2r-1}^{(i_2)}\right)\right)
\Biggl.
\Biggr),
\end{equation}

\vspace{4mm}
\begin{equation}
\label{46}
I_{(2)T,t}^{*(i_1)}=
(T-t)^{5/2}\Biggl(
\frac{1}{3}\zeta_0^{(i_1)}+\frac{1}{\sqrt{2}\pi^2}
\sum_{r=1}^{\infty}\frac{1}{r^2}\zeta_{2r}^{(i_1)}-
\frac{1}{\sqrt{2}\pi}\sum_{r=1}^{\infty}
\frac{1}{r}\zeta_{2r-1}^{(i_1)}\Biggr),
\end{equation}

\vspace{6mm}
\noindent
where
$i_1, i_2=1,\ldots,m;$

\vspace{-2mm}
$$
\zeta_{j}^{(i)}=
\int\limits_t^T \phi_{j}(s) d{\bf f}_s^{(i)}
$$ 

\vspace{3mm}
\noindent
are independent standard Gaussian random variables
for various
$i$ or $j,$ and

\begin{equation}
\label{666.6}
\phi_j(s)=\frac{1}{\sqrt{T-t}}
\begin{cases}
1\ &\hbox{when}\ j=0\cr\cr
\sqrt{2}{\rm sin}(2\pi r(s-t)/(T-t))\ &\hbox{when}\ j=2r-1\cr\cr
\sqrt{2}{\rm cos}(2\pi r(s-t)/(T-t))\ &\hbox{when}\ j=2r
\end{cases},
\end{equation}

\vspace{5mm}
\noindent
where $r=1, 2,\ldots$

Obviously, that $I_{(1)T,t}^{*(i_1)},$ $I_{(2)T,t}^{*(i_1)}$ have 
Gaussian distribution and the 
expansions (\ref{42}), (\ref{46}) are too complex
for such simple stochastic integrals 
as $I_{(1)T,t}^{*(i_1)},$ $I_{(2)T,t}^{*(i_1)}$.

Milstein G.N. proposed \cite{Mi2} the following mean-square 
approximations on the base of the expansions
(\ref{42}), (\ref{43})

\vspace{1mm}
\begin{equation}
\label{444}
I_{(1)T,t}^{*(i_1)q}=-\frac{{(T-t)}^{3/2}}{2}
\Biggl(\zeta_0^{(i_1)}-\frac{\sqrt{2}}{\pi}\biggl(\sum_{r=1}^{q}
\frac{1}{r}
\zeta_{2r-1}^{(i_1)}+\sqrt{\alpha_q}\xi_q^{(i_1)}\biggr)
\Biggr),
\end{equation}

\vspace{4mm}
$$
I_{(00)T,t}^{*(i_1 i_2)q}=\frac{1}{2}(T-t)\Biggl(
\zeta_{0}^{(i_1)}\zeta_{0}^{(i_2)}
+\frac{1}{\pi}
\sum_{r=1}^{q}\frac{1}{r}\left(
\zeta_{2r}^{(i_1)}\zeta_{2r-1}^{(i_2)}-
\zeta_{2r-1}^{(i_1)}\zeta_{2r}^{(i_2)}+
\right.\Biggr.
$$

\vspace{2mm}
\begin{equation}
\label{555}
+\Biggl.\left.\sqrt{2}\left(\zeta_{2r-1}^{(i_1)}\zeta_{0}^{(i_2)}-
\zeta_{0}^{(i_1)}\zeta_{2r-1}^{(i_2)}\right)\right)
+\frac{\sqrt{2}}{\pi}\sqrt{\alpha_q}\left(
\xi_q^{(i_1)}\zeta_0^{(i_2)}-\zeta_0^{(i_1)}\xi_q^{(i_2)}\right)\Biggr),
\end{equation}

\vspace{5mm}
\noindent
where

\vspace{-1mm}
\begin{equation}
\label{333}
\xi_q^{(i)}=\frac{1}{\sqrt{\alpha_q}}\sum_{r=q+1}^{\infty}
\frac{1}{r}~\zeta_{2r-1}^{(i)},\ \ \ 
\alpha_q=\frac{\pi^2}{6}-\sum_{r=1}^q\frac{1}{r^2},
\end{equation}

\vspace{5mm}
\noindent
where
$\zeta_0^{(i)},$ $\zeta_{2r}^{(i)},$
$\zeta_{2r-1}^{(i)},$ $\xi_q^{(i)};$ $r=1,\ldots,q;$
$i=1,\ldots,m$
are independent standard Gaussian random variables.

The approximation $I_{(2)T,t}^{*(i_1)q},$ which corresponds 
to (\ref{444}), (\ref{555})
has the form \cite{KlPl2}

\vspace{1mm}
$$
I_{(2)T,t}^{*(i_1)q}=
(T-t)^{5/2}\Biggl(
\frac{1}{3}\zeta_0^{(i_1)}+\frac{1}{\sqrt{2}\pi^2}\Biggl(
\sum_{r=1}^{q}\frac{1}{r^2}\zeta_{2r}^{(i_1)}+
\sqrt{\beta_q}\mu_q^{(i_1)}\Biggr)-\Biggr.
$$

\vspace{1mm}
\begin{equation}
\label{1970}
\Biggl.-
\frac{1}{\sqrt{2}\pi}\Biggl(\sum_{r=1}^q
\frac{1}{r}\zeta_{2r-1}^{(i_1)}+\sqrt{\alpha_q}\xi_q^{(i_1)}\Biggr)\Biggr),
\end{equation}

\vspace{5mm}
\noindent
where
$\xi_q^{(i)},$ $\alpha_q$ has the form (\ref{333}) and 

\vspace{1mm}
$$
\mu_q^{(i)}=\frac{1}{\sqrt{\beta_q}}\sum_{r=q+1}^{\infty}
\frac{1}{r^2}~\zeta_{2r}^{(i)},\ \ \
\beta_q=\frac{\pi^4}{90}-\sum_{r=1}^q\frac{1}{r^4},
$$

\vspace{5mm}
\noindent
$\phi_j(s)$ is defined by (\ref{666.6}); 
$\zeta_0^{(i)},$ $\zeta_{2r}^{(i)},$
$\zeta_{2r-1}^{(i)},$ $\xi_q^{(i)},$ $\mu_q^{(i)};$ $r=1,\ldots,q;$
$i=1,\ldots,m$ are independent
standard Gaussian random variables;
$i=1,\ldots,m.$

Nevetheless, the expansions (\ref{444}), (\ref{1970}) are too complex for
the approximation of two Gaussian random variables
$I_{(1)T,t}^{*(i_1)},$ $I_{(2)T,t}^{*(i_1)}$.

Using Theorems 1--3 and complete orthonormal system 
of Legendre polynomials
in the space $L_2([t, T])$, we obtain  
for $i_1, i_2=1,\ldots,m$ \cite{2006}-\cite{arxiv-8}

\begin{equation}
\label{4001}
I_{(0)T,t}^{*(i_1)}=\sqrt{T-t}\zeta_0^{(i_1)},
\end{equation}

\vspace{1mm}
\begin{equation}
\label{4002}
I_{(1)T,t}^{*(i_1)}=-\frac{(T-t)^{3/2}}{2}\biggl(\zeta_0^{(i_1)}+
\frac{1}{\sqrt{3}}\zeta_1^{(i_1)}\biggr),
\end{equation}

\vspace{1mm}
\begin{equation}
\label{4003}
I_{(2)T,t}^{*(i_1)}=\frac{(T-t)^{5/2}}{3}\biggl(\zeta_0^{(i_1)}+
\frac{\sqrt{3}}{2}\zeta_1^{(i_1)}+
\frac{1}{2\sqrt{5}}\zeta_2^{(i_1)}\biggr),
\end{equation}

\vspace{2mm}
\begin{equation}
\label{4004}
I_{(00)T,t}^{*(i_1 i_2)}=
\frac{T-t}{2}\Biggl(\zeta_0^{(i_1)}\zeta_0^{(i_2)}+\sum_{i=1}^{\infty}
\frac{1}{\sqrt{4i^2-1}}\left(
\zeta_{i-1}^{(i_1)}\zeta_{i}^{(i_2)}-
\zeta_i^{(i_1)}\zeta_{i-1}^{(i_2)}\right)\Biggr),
\end{equation}

\vspace{4mm}
\noindent
$$
\zeta_{j}^{(i)}=
\int\limits_t^T \phi_{j}(s) d{\bf f}_s^{(i)}
$$ 

\vspace{3mm}
\noindent
are independent standard Gaussian random variables
for various
$i$ or $j$,
where

\vspace{1mm}
\begin{equation}
\label{4009}
\phi_j(x)=\sqrt{\frac{2j+1}{T-t}}P_j\left(\left(
x-\frac{T+t}{2}\right)\frac{2}{T-t}\right);\ j=0, 1, 2,\ldots,
\end{equation}

\vspace{4mm}
\noindent
where $P_j(x)$ is the Legendre polynomial.

It is not difficult to see that the expansions 
(\ref{4002}), (\ref{4003}) are much 
simpler than the expansions (\ref{444}), (\ref{1970}).

Obviously that the 
Milstein approach \cite{Mi2} leads
to iterated series (iterated application of the operation
of limit transitions) 
in contradiction to multiple series 
(the operation of limit transition is implemented only once) 
from Theorems 1--7.

For the case of simplest 
stochastic integral $I_{(00)T,t}^{*(i_1 i_2)}$ of second multiplicity 
this problem was avoided 
as we saw earlier.  However, the situation is not 
the same for the simplest iterated stochastic integral
$I_{(000)T,t}^{*(i_1 i_2 i_3)}$ of third multiplicity.

Let us denote

\vspace{-1mm}
$$
{J}_{(\lambda_{1}\ldots \lambda_k)T,t}^{*(i_1\ldots
i_k)}=
{\int\limits_t^{*}}^T\ldots
{\int\limits_t^{*}}^{t_2}
d{\bf w}_{t_{1}}^{(i_1)}\ldots
d{\bf w}_{t_k}^{(i_k)},
$$

\vspace{4mm}
\noindent
where $\lambda_l=1$ if 
$i_l=1,\ldots,m$ and 
$\lambda_l=0$ if $i_l=0;$ $l=1,\ldots,k$
(${\bf w}_{\tau}^{(i)}={\bf f}_{\tau}^{(i)}$
for $i=1,\ldots,m$ and
${\bf w}_{\tau}^{(0)}=\tau$).

Consider the expansion of iterated Stratonovich stochastic
integral of third multiplicity obtained in \cite{KlPl2}-\cite{Zapad-9}
by the Milstein approach

\vspace{2mm}
$$
J_{(111)\Delta,0}^{*(i_1 i_2 i_3)}
=\frac{1}{\Delta}
J_{(1)\Delta,0}^{*(i_1)}J_{(011)\Delta,0}^{*(0 i_2 i_3)}
+\frac{1}{2}a_{i_1,0}J_{(11)\Delta,0}^{*(i_2 i_3)}+
\frac{1}{2\pi}b_{i_1}J_{(1)\Delta,0}^{(i_2)}
J_{(1)\Delta,0}^{*(i_3)}-
$$

\vspace{3mm}
\begin{equation}
\label{starr}
-\Delta J_{(1)\Delta,0}^{*(i_2)}
B_{i_1 i_3}+\Delta J_{(1)\Delta,0}^{*(i_3)}
\left(\frac{1}{2}A_{i_1i_2}-C_{i_2i_1}\right)
+\Delta^{3/2}D_{i_1i_2i_3},
\end{equation}

\vspace{6mm}
\noindent
where

\vspace{-3mm}
$$
J_{(011)\Delta,0}^{*(0 i_2 i_3)}=\frac{1}{6}
J_{(1)\Delta,0}^{*(i_2)}
J_{(1)\Delta,0}^{*(i_3)}-\frac{1}{\pi}\Delta
J_{(1)\Delta,0}^{*(i_3)}b_{i_2}+
$$

\vspace{3mm}
$$
+\Delta^2 B_{i_2i_3}-\frac{1}{4}\Delta a_{i_3,0}
J_{(1)\Delta,0}^{*(i_2)}+
\frac{1}{2\pi}\Delta b_{i_3}
J_{(1)\Delta,0}^{*(i_2)}
+\Delta^2 C_{i_2i_3}+\frac{1}{2}\Delta^2 A_{i_2i_3},
$$

\vspace{6mm}

$$
A_{i_2i_3}=\frac{\pi}{\Delta}\sum_{r=1}^{\infty}
r\left(a_{i_2,r}b_{i_3,r}-
b_{i_2,r}
a_{i_3,r}\right),
$$

\vspace{4mm}

$$
C_{i_2i_3}=-\frac{1}{\Delta}\sum_{l=1}^{\infty}
\sum_{r=1 (r\ne l)}^{\infty}
\frac{r}{r^2-l^2}\left(ra_{i_2,r}a_{i_3,l}+
lb_{i_2,r}
b_{i_3,l}\right),
$$

\vspace{4mm}

$$
B_{i_2i_3}=\frac{1}{2\Delta}\sum_{r=1}^{\infty}
\left(a_{i_2,r}a_{i_3,r}+
b_{i_2,r}
b_{i_3,r}\right),\
b_i=\sum_{r=1}^{\infty}\frac{1}{r}b_{i,r},
$$

\vspace{9mm}

$$
D_{i_1i_2i_3}=-\frac{\pi}{2\Delta^{3/2}}
\sum_{l=1}^{\infty}\sum_{r=1}^{\infty}
l
\Biggl(a_{i_2,l}
\left(a_{i_3,l+r}b_{i_1,r}-a_{i_1,r}
b_{i_3,l+r}\right)
+\Biggr.
$$

\vspace{1mm}
$$
\Biggl.
+
b_{i_2,l}
\left(a_{i_1,r}a_{i_3,r+l}+b_{i_1,r}
b_{i_3,l+r}\right)\Biggr)+
$$

\vspace{1mm}
$$
+\frac{\pi}{2\Delta^{3/2}}
\sum_{l=1}^{\infty}\sum_{r=1}^{l-1}
l
\Biggl(a_{i_2,l}
\left(a_{i_1,r}b_{i_3,l-r}+a_{i_3,l-r}
b_{i_1,r}\right)
-\Biggr.
$$

\vspace{1mm}
$$
\Biggl.
-b_{i_2,l}
\left(a_{i_1,r}a_{i_3,l-r}-b_{i_1,r}
b_{i_3,l-r}\right)\Biggr)+
$$

\vspace{1mm}
$$
+\frac{\pi}{2\Delta^{3/2}}
\sum_{l=1}^{\infty}\sum_{r=l+1}^{\infty}
l
\Biggl(a_{i_2,l}
\left(a_{i_3,r-l}b_{i_1,r}-a_{i_1,r}
b_{i_3,r-l}\right)
+\Biggr.
$$

\vspace{1mm}
$$
\Biggl.
+
b_{i_2,l}
\left(a_{i_1,r}a_{i_3,r-l}+b_{i_1,r}
b_{i_3,r-l}\right)\Biggr).
$$

\vspace{5mm}

From 
the form
of expansion
(\ref{starr}) and expansion 
of the stochastic integral
$J_{(011)\Delta,0}^{*(0i_2 i_3)}$
we can conclude
that they include 
iterated (double)
series. 
Moreover, for approximation of the considered
stochastic integral $J_{(111)\Delta,0}^{*(i_1 i_2 i_3)}$ in the works
\cite{KlPl2}
(Sect.~5.8, pp.~202--204), \cite{KPS} (pp.~82-84),
\cite{KPW} (pp.~438-439),  
\cite{Zapad-9} (pp.~263-264)
it is proposed to put upper limits of summation 
by equal $q$
(on the base of the Wong--Zakai approximation \cite{W-Z-1}-\cite{Watanabe},
but without rigorous proof;
also see discussion in Sect.~7).

For example, the value $D_{i_1i_2i_3}$ is approximated 
in \cite{KlPl2}
(Sect.~5.8, pp.~202--204), \cite{KPS} (pp.~82-84),
\cite{KPW} (pp.~438-439),  
\cite{Zapad-9} (pp.~263-264)
by the double sums of the form

\vspace{2mm}
$$
D_{i_1i_2i_3}^{(q)}=-\frac{\pi}{2\Delta^{3/2}}
\sum_{l=1}^{q}\sum_{r=1}^{q}
l
\Biggl(a_{i_2,l}
\left(a_{i_3,l+r}b_{i_1,r}-a_{i_1,r}
b_{i_3,l+r}\right)
+\Biggr.
$$

\vspace{1mm}
$$
\Biggl.
+
b_{i_2,l}
\left(a_{i_1,r}a_{i_3,r+l}+b_{i_1,r}
b_{i_3,l+r}\right)\Biggr)+
$$

\vspace{1mm}
$$
+\frac{\pi}{2\Delta^{3/2}}
\sum_{l=1}^{q}\sum_{r=1}^{l-1}
l
\Biggl(a_{i_2,l}
\left(a_{i_1,r}b_{i_3,l-r}+a_{i_3,l-r}
b_{i_1,r}\right)
-\Biggr.
$$

\vspace{1mm}
$$
\Biggl.
-b_{i_2,l}
\left(a_{i_1,r}a_{i_3,l-r}-b_{i_1,r}
b_{i_3,l-r}\right)\Biggr)+
$$

\vspace{1mm}
$$
+\frac{\pi}{2\Delta^{3/2}}
\sum_{l=1}^{q}\sum_{r=l+1}^{2q}
l
\Biggl(a_{i_2,l}
\left(a_{i_3,r-l}b_{i_1,r}-a_{i_1,r}
b_{i_3,r-l}\right)
+\Biggr.
$$

\vspace{1mm}
$$
\Biggl.
+
b_{i_2,l}
\left(a_{i_1,r}a_{i_3,r-l}+b_{i_1,r}
b_{i_3,r-l}\right)\Biggr).
$$

\vspace{5mm}

Obviously, we can avoid 
this problem (iterated application of the operation
of limit transition)  
using the method based on Theorems 1--7.

If we prove that the terms of the expansion (\ref{starr}) coincide 
with the terms of its analogue obtained using Theorems 1--3 
(this fact is proved in \cite{2006}-\cite{2018axx} for the 
simplest stochastic integrals $I_{(1)T,t}^{*(i_1)},$ 
$I_{(00)T,t}^{*(i_1i_2)}$
of first and second multiplicity), 
then we can replace the iterated (double) series in (\ref{starr}) by 
the multiple 
ones, as in Theorems 1--3 (as was made formally 
in \cite{KlPl2}-\cite{Zapad-9}). 
However, it 
requires a separate argumentation.

\vspace{5mm}

\section{Approximation of Specific Iterated Stochastic Integrals of 
Multiplicities 1 to 3 Using Theorem 3 and 
Trigonometric System of Functions}

\vspace{5mm}

In \cite{2006}-\cite{2018axx}
on the base of Theorems 1--3 
the author of this paper obtained
the
following expansions of the iterated Stratonovich stochastic
integrals (\ref{k1001}) 
(independently from the papers
\cite{Mi2}-\cite{Mi3}
excepting the method in which additional
random variables $\xi_q^{(i)}$ and $\mu_q^{(i)})$ 
are introduced)

\vspace{1mm}
 
\begin{equation}
\label{410}
I_{(0)T,t}^{*(i_1)}=\sqrt{T-t}\zeta_0^{(i_1)},
\end{equation}

\vspace{1mm}
\begin{equation}
\label{420}
I_{(1)T,t}^{*(i_1)q}=-\frac{{(T-t)}^{3/2}}{2}
\Biggl(\zeta_0^{(i_1)}-\frac{\sqrt{2}}{\pi}\Biggl(\sum_{r=1}^{q}
\frac{1}{r}
\zeta_{2r-1}^{(i_1)}+\sqrt{\alpha_q}\xi_q^{(i_1)}\Biggr)
\Biggr),
\end{equation}

\vspace{4mm}

$$
I_{(00)T,t}^{*(i_1 i_2)q}=\frac{1}{2}(T-t)\Biggl(
\zeta_{0}^{(i_1)}\zeta_{0}^{(i_2)}
+\frac{1}{\pi}
\sum_{r=1}^{q}\frac{1}{r}\left(
\zeta_{2r}^{(i_1)}\zeta_{2r-1}^{(i_2)}-
\zeta_{2r-1}^{(i_1)}\zeta_{2r}^{(i_2)}+
\right.\Biggr.
$$

\vspace{1mm}
\begin{equation}
\label{430}
+\Biggl.\left.\sqrt{2}\left(\zeta_{2r-1}^{(i_1)}\zeta_{0}^{(i_2)}-
\zeta_{0}^{(i_1)}\zeta_{2r-1}^{(i_2)}\right)\right)
+\frac{\sqrt{2}}{\pi}\sqrt{\alpha_q}\left(
\xi_q^{(i_1)}\zeta_0^{(i_2)}-\zeta_0^{(i_1)}\xi_q^{(i_2)}\right)\Biggr),
\end{equation}

\vspace{6mm}

$$
I_{(000)T,t}^{*(i_1 i_2 i_3)q}=(T-t)^{3/2}\Biggl(\frac{1}{6}
\zeta_{0}^{(i_1)}\zeta_{0}^{(i_2)}\zeta_{0}^{(i_3)}+\Biggr.
\frac{\sqrt{\alpha_q}}{2\sqrt{2}\pi}\left(
\xi_q^{(i_1)}\zeta_0^{(i_2)}\zeta_0^{(i_3)}-\xi_q^{(i_3)}\zeta_0^{(i_2)}
\zeta_0^{(i_1)}\right)+
$$
$$
+\frac{1}{2\sqrt{2}\pi^2}\sqrt{\beta_q}\left(
\mu_q^{(i_1)}\zeta_0^{(i_2)}\zeta_0^{(i_3)}-2\mu_q^{(i_2)}\zeta_0^{(i_1)}
\zeta_0^{(i_3)}+\mu_q^{(i_3)}\zeta_0^{(i_1)}\zeta_0^{(i_2)}\right)+
$$
$$
+
\frac{1}{2\sqrt{2}}\sum_{r=1}^{q}
\Biggl(\frac{1}{\pi r}\left(
\zeta_{2r-1}^{(i_1)}
\zeta_{0}^{(i_2)}\zeta_{0}^{(i_3)}-
\zeta_{2r-1}^{(i_3)}
\zeta_{0}^{(i_2)}\zeta_{0}^{(i_1)}\right)+\Biggr.
$$
$$
\Biggl.+
\frac{1}{\pi^2 r^2}\left(
\zeta_{2r}^{(i_1)}
\zeta_{0}^{(i_2)}\zeta_{0}^{(i_3)}-
2\zeta_{2r}^{(i_2)}
\zeta_{0}^{(i_3)}\zeta_{0}^{(i_1)}+
\zeta_{2r}^{(i_3)}
\zeta_{0}^{(i_2)}\zeta_{0}^{(i_1)}\right)\Biggr)+
$$
$$
+
\sum_{r=1}^{q}
\Biggl(\frac{1}{4\pi r}\left(
\zeta_{2r}^{(i_1)}
\zeta_{2r-1}^{(i_2)}\zeta_{0}^{(i_3)}-
\zeta_{2r-1}^{(i_1)}
\zeta_{2r}^{(i_2)}\zeta_{0}^{(i_3)}-
\zeta_{2r-1}^{(i_2)}
\zeta_{2r}^{(i_3)}\zeta_{0}^{(i_1)}+
\zeta_{2r-1}^{(i_3)}
\zeta_{2r}^{(i_2)}\zeta_{0}^{(i_1)}\right)+\Biggr.
$$
$$
+
\frac{1}{8\pi^2 r^2}\left(
3\zeta_{2r-1}^{(i_1)}
\zeta_{2r-1}^{(i_2)}\zeta_{0}^{(i_3)}+
\zeta_{2r}^{(i_1)}
\zeta_{2r}^{(i_2)}\zeta_{0}^{(i_3)}-
6\zeta_{2r-1}^{(i_1)}
\zeta_{2r-1}^{(i_3)}\zeta_{0}^{(i_2)}+\right.
$$
\begin{equation}
\label{44}
\Biggl.\left.
+
3\zeta_{2r-1}^{(i_2)}
\zeta_{2r-1}^{(i_3)}\zeta_{0}^{(i_1)}-
2\zeta_{2r}^{(i_1)}
\zeta_{2r}^{(i_3)}\zeta_{0}^{(i_2)}+
\zeta_{2r}^{(i_3)}
\zeta_{2r}^{(i_2)}\zeta_{0}^{(i_1)}\right)\Biggr)
\Biggl.+D_{T,t}^{(i_1i_2i_3)q}\Biggr),
\end{equation}

\vspace{6mm}
\noindent
where

$$
D_{T,t}^{(i_1i_2i_3)q}=
\frac{1}{2\pi^2}\sum_{\stackrel{r,l=1}{{}_{r\ne l}}}^{q}
\Biggl(\frac{1}{r^2-l^2}\biggl(
\zeta_{2r}^{(i_1)}
\zeta_{2l}^{(i_2)}\zeta_{0}^{(i_3)}-
\zeta_{2r}^{(i_2)}
\zeta_{0}^{(i_1)}\zeta_{2l}^{(i_3)}+\biggr.\Biggr.
$$
$$
\Biggl.+\biggl.
\frac{r}{l}
\zeta_{2r-1}^{(i_1)}
\zeta_{2l-1}^{(i_2)}\zeta_{0}^{(i_3)}-\frac{l}{r}
\zeta_{0}^{(i_1)}
\zeta_{2r-1}^{(i_2)}\zeta_{2l-1}^{(i_3)}\biggr)-
\frac{1}{rl}\zeta_{2r-1}^{(i_1)}
\zeta_{0}^{(i_2)}\zeta_{2l-1}^{(i_3)}\Biggr)+
$$
$$
+
\frac{1}{4\sqrt{2}\pi^2}\Biggl(
\sum_{r,m=1}^{q}\Biggl(\frac{2}{rm}
\left(-\zeta_{2r-1}^{(i_1)}
\zeta_{2m-1}^{(i_2)}\zeta_{2m}^{(i_3)}+
\zeta_{2r-1}^{(i_1)}
\zeta_{2r}^{(i_2)}\zeta_{2m-1}^{(i_3)}+
\right.\Biggr.\Biggr.
$$

$$
\left.+
\zeta_{2r-1}^{(i_1)}
\zeta_{2m}^{(i_2)}\zeta_{2m-1}^{(i_3)}-
\zeta_{2r}^{(i_1)}
\zeta_{2r-1}^{(i_2)}\zeta_{2m-1}^{(i_3)}\right)+
$$

$$
+\frac{1}{m(r+m)}
\left(-\zeta_{2(m+r)}^{(i_1)}
\zeta_{2r}^{(i_2)}\zeta_{2m}^{(i_3)}-
\zeta_{2(m+r)-1}^{(i_1)}
\zeta_{2r-1}^{(i_2)}\zeta_{2m}^{(i_3)}-
\right.
$$
$$
\Biggl.\left.
-\zeta_{2(m+r)-1}^{(i_1)}
\zeta_{2r}^{(i_2)}\zeta_{2m-1}^{(i_3)}+
\zeta_{2(m+r)}^{(i_1)}
\zeta_{2r-1}^{(i_2)}\zeta_{2m-1}^{(i_3)}\right)\Biggr)+
$$
$$
+
\sum_{m=1}^{q}\sum_{l=m+1}^{q}\Biggl(\frac{1}{m(l-m)}
\left(\zeta_{2(l-m)}^{(i_1)}
\zeta_{2l}^{(i_2)}\zeta_{2m}^{(i_3)}+
\zeta_{2(l-m)-1}^{(i_1)}
\zeta_{2l-1}^{(i_2)}\zeta_{2m}^{(i_3)}-
\right.\Biggr.
$$

$$
\left.
-\zeta_{2(l-m)-1}^{(i_1)}
\zeta_{2l}^{(i_2)}\zeta_{2m-1}^{(i_3)}+
\zeta_{2(l-m)}^{(i_1)}
\zeta_{2l-1}^{(i_2)}\zeta_{2m-1}^{(i_3)}\right)+
$$

$$
+
\frac{1}{l(l-m)}
\left(-\zeta_{2(l-m)}^{(i_1)}
\zeta_{2m}^{(i_2)}\zeta_{2l}^{(i_3)}+
\zeta_{2(l-m)-1}^{(i_1)}
\zeta_{2m-1}^{(i_2)}\zeta_{2l}^{(i_3)}-
\right.
$$
$$
\Biggl.
\Biggl.
\Biggl.
\left.
-\zeta_{2(l-m)-1}^{(i_1)}
\zeta_{2m}^{(i_2)}\zeta_{2l-1}^{(i_3)}-
\zeta_{2(l-m)}^{(i_1)}
\zeta_{2m-1}^{(i_2)}\zeta_{2l-1}^{(i_3)}\right)\Biggr)\Biggr),
$$

\vspace{10mm}

$$
I_{(10)T,t}^{*(i_1 i_2)q}=-(T-t)^{2}\Biggl(\frac{1}{6}
\zeta_{0}^{(i_1)}\zeta_{0}^{(i_2)}-\frac{1}{2\sqrt{2}\pi}
\sqrt{\alpha_q}\xi_q^{(i_2)}\zeta_0^{(i_1)}+\Biggr.
$$
$$
+\frac{1}{2\sqrt{2}\pi^2}\sqrt{\beta_q}\Biggl(
\mu_q^{(i_2)}\zeta_0^{(i_1)}-2\mu_q^{(i_1)}\zeta_0^{(i_2)}\Biggr)+
$$
$$
+\frac{1}{2\sqrt{2}}\sum_{r=1}^{q}
\Biggl(-\frac{1}{\pi r}
\zeta_{2r-1}^{(i_2)}
\zeta_{0}^{(i_1)}+
\frac{1}{\pi^2 r^2}\left(
\zeta_{2r}^{(i_2)}
\zeta_{0}^{(i_1)}-
2\zeta_{2r}^{(i_1)}
\zeta_{0}^{(i_2)}\right)\Biggr)-
$$
$$
-
\frac{1}{2\pi^2}\sum_{\stackrel{r,l=1}{{}_{r\ne l}}}^{q}
\frac{1}{r^2-l^2}\Biggl(
\zeta_{2r}^{(i_1)}
\zeta_{2l}^{(i_2)}+
\frac{l}{r}
\zeta_{2r-1}^{(i_1)}
\zeta_{2l-1}^{(i_2)}
\Biggr)+
$$
\begin{equation}
\label{9440}
+
\sum_{r=1}^{q}
\Biggl(\frac{1}{4\pi r}\left(
\zeta_{2r}^{(i_1)}
\zeta_{2r-1}^{(i_2)}-
\zeta_{2r-1}^{(i_1)}
\zeta_{2r}^{(i_2)}\right)+
\Biggl.\Biggl.
\frac{1}{8\pi^2 r^2}\left(
3\zeta_{2r-1}^{(i_1)}
\zeta_{2r-1}^{(i_2)}+
\zeta_{2r}^{(i_2)}
\zeta_{2r}^{(i_1)}\right)\Biggr)\Biggr),
\end{equation}

\vspace{8mm}

$$
I_{(01)T,t}^{*(i_1 i_2)q}=
(T-t)^{2}\Biggl(-\frac{1}{3}
\zeta_{0}^{(i_1)}\zeta_{0}^{(i_2)}-\frac{1}{2\sqrt{2}\pi}
\sqrt{\alpha_q}\left(\xi_q^{(i_1)}\zeta_0^{(i_2)}-
2\xi_q^{(i_2)}\zeta_0^{(i_1)}\right)
+\Biggr.
$$
$$
+\frac{1}{2\sqrt{2}\pi^2}\sqrt{\beta_q}\Biggl(
\mu_q^{(i_1)}\zeta_0^{(i_2)}-2\mu_q^{(i_2)}\zeta_0^{(i_1)}\Biggr)-
$$
$$
-\frac{1}{2\sqrt{2}}\sum_{r=1}^{q}
\Biggl(\frac{1}{\pi r}\left(
\zeta_{2r-1}^{(i_1)}
\zeta_{0}^{(i_2)}-
2\zeta_{2r-1}^{(i_2)}
\zeta_{0}^{(i_1)}\right)
-\frac{1}{\pi^2 r^2}\left(
\zeta_{2r}^{(i_1)}
\zeta_{0}^{(i_2)}-
2\zeta_{2r}^{(i_2)}
\zeta_{0}^{(i_1)}\right)\Biggr)+
$$
$$
+
\frac{1}{2\pi^2}\sum_{\stackrel{r,l=1}{{}_{r\ne l}}}^{q}
\frac{1}{r^2-l^2}\Biggl(
\frac{r}{l}\zeta_{2r-1}^{(i_1)}
\zeta_{2l-1}^{(i_2)}+
\zeta_{2r}^{(i_1)}
\zeta_{2l}^{(i_2)}
\Biggr)-
$$
\begin{equation}
\label{450}
-
\sum_{r=1}^{q}
\Biggl(\frac{1}{4\pi r}\left(
\zeta_{2r}^{(i_1)}
\zeta_{2r-1}^{(i_2)}-
\zeta_{2r-1}^{(i_1)}
\zeta_{2r}^{(i_2)}\right)-
\Biggl.\Biggl.
\frac{1}{8\pi^2 r^2}\left(
3\zeta_{2r-1}^{(i_1)}
\zeta_{2r-1}^{(i_2)}+
\zeta_{2r}^{(i_1)}
\zeta_{2r}^{(i_2)}\right)\Biggr)\Biggr),
\vspace{8mm}
\end{equation}

\vspace{1mm}

$$
I_{(2)T,t}^{*(i_1)q}=
(T-t)^{5/2}\Biggl(
\frac{1}{3}\zeta_0^{(i_1)}+\frac{1}{\sqrt{2}\pi^2}\Biggl(
\sum_{r=1}^{q}\frac{1}{r^2}\zeta_{2r}^{(i_1)}+
\sqrt{\beta_q}\mu_q^{(i_1)}\Biggr)-\Biggr.
$$
\begin{equation}
\label{460}
\Biggl.-
\frac{1}{\sqrt{2}\pi}\Biggl(\sum_{r=1}^q
\frac{1}{r}\zeta_{2r-1}^{(i_1)}+\sqrt{\alpha_q}\xi_q^{(i_1)}\Biggr)\Biggr),
\end{equation}

\vspace{4mm}
\noindent
where

$$
\xi_q^{(i)}=\frac{1}{\sqrt{\alpha_q}}\sum_{r=q+1}^{\infty}
\frac{1}{r}~\zeta_{2r-1}^{(i)},\ \ \
\alpha_q=\frac{\pi^2}{6}-\sum_{r=1}^q\frac{1}{r^2},\ \ \
\mu_q^{(i)}=\frac{1}{\sqrt{\beta_q}}\sum_{r=q+1}^{\infty}
\frac{1}{r^2}~\zeta_{2r}^{(i)},
$$

$$
\beta_q=\frac{\pi^4}{90}-\sum_{r=1}^q\frac{1}{r^4},\ \ \
\zeta_j^{(i)}=\int\limits_t^T\phi_j(s)d{\bf f}_s^{(i)},
$$

\vspace{7mm}
\noindent
where
$\phi_j(s)$ has the form (\ref{666.6}); 
$\zeta_0^{(i)},$ $\zeta_{2r}^{(i)},$
$\zeta_{2r-1}^{(i)},$ $\xi_q^{(i)},$ $\mu_q^{(i)};$ $r=1,\ldots,q;$
$i=1,\ldots,m$ are independent
standard Gaussian random variables;
$i_1, i_2, i_3=1,\ldots,m.$

Note that from (\ref{9440}), (\ref{450}) it follows that

\begin{equation}
\label{leto3000mil}
\sum\limits_{j=0}^{\infty}C_{jj}^{10}=
\sum\limits_{j=0}^{\infty}C_{jj}^{01}=-\frac{(T-t)^2}{4},
\end{equation}

\vspace{2mm}
\noindent
where

$$
C_{jj}^{10}=\int\limits_{t}^{T}\phi_j(x)
\int\limits_{t}^{x}\phi_j(y)(t-y)
dy dx,
$$

$$
C_{jj}^{01}=\int\limits_{t}^{T}\phi_j(x)(t-x)
\int\limits_{t}^{x}\phi_j(y)
dy dx.
$$

\vspace{5mm}

The formulas (\ref{leto3000mil}) are particular cases 
of the more general relation, which we applied for the proof 
of Theorem 3 for the case $k=2$ (see \cite{2010-2}-\cite{2018axx}).

Let us consider the mean-square errors of approximations
(\ref{430})--(\ref{450}). From the relations (\ref{430})--(\ref{450})
when $i_1\ne i_2,$ $i_2\ne i_3,$ $i_1\ne i_3$ we obtain
by direct 
calculation 

\vspace{1mm}
\begin{equation}
\label{801}
{\sf M}\left\{\left(I_{(00)T,t}^{*(i_1 i_2)}-
I_{(00)T,t}^{*(i_1 i_2)q}
\right)^2\right\}
=\frac{(T-t)^{2}}{2\pi^2}\Biggl(\frac{\pi^2}{6}-
\sum_{r=1}^q \frac{1}{r^2}\Biggr),
\end{equation}

\vspace{5mm}

$$
{\sf M}\left\{\left(I_{(000)T,t}^{*(i_1 i_2 i_3)}-
I_{(000)T,t}^{*(i_1 i_2 i_3)q}\right)^2\right\}
=(T-t)^{3}\Biggl(\frac{1}{4\pi^2}
\Biggl(\frac{\pi^2}{6}-
\sum_{r=1}^q \frac{1}{r^2}\Biggr)+
\Biggr.
$$

\vspace{2mm}
\begin{equation}
\label{802}
\Biggl.
+\frac{55}{32\pi^4}\Biggl(\frac{\pi^4}{90}-
\sum_{r=1}^q \frac{1}{r^4}\Biggr)
+\frac{1}{4\pi^4}
\Biggl(\sum_{\stackrel{r,l=1}{{}_{r\ne l}}}^{\infty}
-\sum_{\stackrel{r,l=1}{{}_{r\ne l}}}^{q}
\Biggr)
\frac{5l^4+4r^4-3l^2r^2}{r^2 l^2(r^2-l^2)^2}\Biggr),
\end{equation}

\vspace{6mm}

$$
{\sf M}\left\{\left(I_{(01)T,t}
^{*(i_1i_2)}-I_{(01)T,t}^{*(i_1i_2)q}\right)^2\right\}
=(T-t)^{4}\Biggl(\frac{1}{8\pi^2}
\Biggl(\frac{\pi^2}{6}-
\sum_{r=1}^q \frac{1}{r^2}\Biggr)+\Biggr.
$$

\vspace{2mm}
\begin{equation}
\label{804}
\Biggl.+\frac{5}{32\pi^4}\Biggl(\frac{\pi^4}{90}-
\sum_{r=1}^q \frac{1}{r^4}\Biggr)+
\frac{1}{4\pi^4}\Biggl(\sum_{\stackrel{k,l=1}{{}_{k\ne l}}}^{\infty}
-\sum_{\stackrel{k,l=1}{{}_{k\ne l}}}^{q}
\Biggr)\frac{l^2+k^2}{k^2(l^2-k^2)^2}\Biggr),
\end{equation}

\vspace{6mm}

$$
{\sf M}\left\{\left(I_{(10)T,t}^{*(i_1i_2)}
-I_{(10)T,t}^{*(i_1i_2)q}\right)^2\right\}
=(T-t)^{4}\Biggl(\frac{1}{8\pi^2}
\Biggl(\frac{\pi^2}{6}-
\sum_{r=1}^q \frac{1}{r^2}\Biggr)+\Biggr.
$$

\vspace{2mm}
\begin{equation}
\label{805}
+\Biggl.
\frac{5}{32\pi^4}\Biggl(\frac{\pi^4}{90}-
\sum_{r=1}^q \frac{1}{r^4}\Biggr)+
\frac{1}{4\pi^4}\Biggl(\sum_{\stackrel{k,l=1}{{}_{k\ne l}}}^{\infty}
-\sum_{\stackrel{k,l=1}{{}_{k\ne l}}}^{q}
\Biggr)\frac{l^2+k^2}{l^2(l^2-k^2)^2}\Biggr).
\end{equation}

\vspace{8mm}

It is easy to demonstrate that the relations 
(\ref{802}), (\ref{804}), and (\ref{805})
can be represented using Theorem 8
in the following form

\vspace{2mm}
$$
{\sf M}\left\{\left(I_{(000)T,t}^{*(i_1 i_2 i_3)}-
I_{(000)T,t}^{*(i_1 i_2 i_3)q}\right)^2\right\}=
(T-t)^3\Biggl(\frac{4}{45}-\frac{1}{4\pi^2}\sum_{r=1}^q\frac{1}{r^2}-
\Biggl.
$$

\vspace{1mm}
\begin{equation}
\label{101.100}
\Biggl.-\frac{55}{32\pi^4}\sum_{r=1}^q\frac{1}{r^4}-
\frac{1}{4\pi^4}\sum_{\stackrel{r,l=1}{{}_{r\ne l}}}^q
\frac{5l^4+4r^4-3r^2l^2}{r^2 l^2 \left(r^2-l^2\right)^2}\Biggr),
\end{equation}

\vspace{5mm}

$$
{\sf M}\left\{\left(I_{(10)T,t}^{*(i_1 i_2)}-
I_{(10)T,t}^{*(i_1 i_2)q}\right)^2\right\}=
\frac{(T-t)^4}{4}\Biggl(\frac{1}{9}-
\frac{1}{2\pi^2}\sum_{r=1}^q \frac{1}{r^2}-\Biggr.
$$

\vspace{1mm}
\begin{equation}
\label{101.101}
\Biggl.-\frac{5}{8\pi^4}\sum_{r=1}^q \frac{1}{r^4}-
\frac{1}{\pi^4}\sum_{\stackrel{k,l=1}{{}_{k\ne l}}}^q
\frac{k^2+l^2}{l^2\left(l^2-k^2\right)^2}\Biggr),
\end{equation}

\vspace{5mm}

$$
{\sf M}\left\{\left(I_{(01)T,t}^{*(i_1 i_2)}-
I_{(01)T,t}^{*(i_1 i_2)q}\right)^2\right\}=
\frac{(T-t)^4}{4}\Biggl(\frac{1}{9}-
\frac{1}{2\pi^2}\sum_{r=1}^q \frac{1}{r^2}-\Biggr.
$$

\vspace{1mm}
\begin{equation}
\label{101.102}
\Biggl.-\frac{5}{8\pi^4}\sum_{r=1}^q \frac{1}{r^4}-
\frac{1}{\pi^4}\sum_{\stackrel{k,l=1}{{}_{k\ne l}}}^q
\frac{l^2+k^2}{k^2\left(l^2-k^2\right)^2}\Biggr).
\end{equation}

\vspace{5mm}

Comparing (\ref{101.100})--(\ref{101.102}) and
(\ref{802})--(\ref{805}), we obtain

\vspace{1mm}
\begin{equation}
\label{101.103}
\sum_{\stackrel{k,l=1}{{}_{k\ne l}}}^{\infty}\frac{l^2+k^2}
{k^2\left(l^2-k^2\right)^2}=
\sum_{\stackrel{k,l=1}{{}_{k\ne l}}}^{\infty}\frac{l^2+k^2}
{l^2\left(l^2-k^2\right)^2}=\frac{\pi^4}{48},
\end{equation}

\vspace{2mm}
\begin{equation}
\label{daug1}
\sum_{\stackrel{r,l=1}{{}_{r\ne l}}}^{\infty}
\frac{5l^4+4r^4-3r^2 l^2}{r^2 l^2\left(r^2-l^2\right)^2}=
\frac{9\pi^4}{80}.
\end{equation}

\vspace{5mm}

Let us consider approximations of the stochastic 
integrals $I_{(10)T,t}^{*(i_1i_1)},$
$I_{(01)T,t}^{*(i_1i_1)}$ and conditions for selecting 
the number $q$ using the trigonometric system of functions

\vspace{3mm}
$$
I_{(10)T,t}^{*(i_1 i_1)q}=-(T-t)^{2}\Biggl(\frac{1}{6}
\left(\zeta_{0}^{(i_1)}\right)^2-\frac{1}{2\sqrt{2}\pi}
\sqrt{\alpha_q}\xi_q^{(i_1)}\zeta_0^{(i_1)}-\Biggr.
$$

$$
-\frac{1}{2\sqrt{2}\pi^2}\sqrt{\beta_q}
\mu_q^{(i_1)}\zeta_0^{(i_1)}
-\frac{1}{2\sqrt{2}}\sum_{r=1}^{q}
\Biggl(\frac{1}{\pi r}
\zeta_{2r-1}^{(i_1)}
\zeta_{0}^{(i_1)}+
\frac{1}{\pi^2 r^2}
\zeta_{2r}^{(i_1)}
\zeta_{0}^{(i_1)}\Biggr)-
$$

$$
-
\frac{1}{2\pi^2}\sum_{\stackrel{r,l=1}{{}_{r\ne l}}}^{q}
\frac{1}{r^2-l^2}\Biggl(
\zeta_{2r}^{(i_1)}
\zeta_{2l}^{(i_1)}+
\frac{l}{r}
\zeta_{2r-1}^{(i_1)}
\zeta_{2l-1}^{(i_1)}
\Biggr)+
$$

$$
\Biggl.+
\frac{1}{8\pi^2}\sum_{r=1}^{q}
\frac{1}{r^2}\left(
3\left(\zeta_{2r-1}^{(i_1)}\right)^2
+
\left(\zeta_{2r}^{(i_1)}\right)^2\right)\Biggr),
$$

\vspace{9mm}

$$
I_{(01)T,t}^{*(i_1 i_1)q}=(T-t)^{2}\Biggl(-\frac{1}{3}
\left(\zeta_{0}^{(i_1)}\right)^2+\frac{1}{2\sqrt{2}\pi}
\sqrt{\alpha_q}\xi_q^{(i_1)}\zeta_0^{(i_1)}-\Biggr.
$$

\vspace{1mm}
$$
-\frac{1}{2\sqrt{2}\pi^2}\sqrt{\beta_q}
\mu_q^{(i_1)}\zeta_0^{(i_1)}
+\frac{1}{2\sqrt{2}}\sum_{r=1}^{q}
\Biggl(\frac{1}{\pi r}
\zeta_{2r-1}^{(i_1)}
\zeta_{0}^{(i_1)}-
\frac{1}{\pi^2 r^2}
\zeta_{2r}^{(i_1)}
\zeta_{0}^{(i_1)}\Biggr)+
$$

\vspace{1mm}
$$
+
\frac{1}{2\pi^2}\sum_{\stackrel{r,l=1}{{}_{r\ne l}}}^{q}
\frac{1}{r^2-l^2}\Biggl(
\zeta_{2r}^{(i_1)}
\zeta_{2l}^{(i_1)}+
\frac{r}{l}
\zeta_{2r-1}^{(i_1)}
\zeta_{2l-1}^{(i_1)}
\Biggr)+
$$

\vspace{1mm}
$$
\Biggl.+
\frac{1}{8\pi^2}\sum_{r=1}^{q}
\frac{1}{r^2}\left(
3\left(\zeta_{2r-1}^{(i_1)}\right)^2
+
\left(\zeta_{2r}^{(i_1)}\right)^2\right)\Biggr).
$$

\vspace{5mm}

Then, we obtain

\vspace{3mm}

$$
{\sf M}\left\{\left(I_{(01)T,t}
^{*(i_1i_1)}-I_{(01)T,t}^{*(i_1i_1)q}\right)^2\right\}=
{\sf M}\left\{\left(I_{(10)T,t}^{*(i_1i_1)}
-I_{(10)T,t}^{*(i_1i_1)q}\right)^2\right\}=
$$

\vspace{2mm}
$$
=\frac{(T-t)^{4}}{4}\Biggl(\frac{2}{\pi^4}\Biggl(\frac{\pi^4}{90}-
\sum_{r=1}^q \frac{1}{r^4}\Biggr)
+\frac{1}{\pi^4}
\Biggl(\frac{\pi^2}{6}-
\sum_{r=1}^q \frac{1}{r^2}\Biggr)^2+\Biggr.
$$

\vspace{2mm}
\begin{equation}
\label{101.104}
\Biggl.+
\frac{1}{\pi^4}\Biggl(\sum_{\stackrel{k,l=1}{{}_{k\ne l}}}^{\infty}
-\sum_{\stackrel{k,l=1}{{}_{k\ne l}}}^{q}
\Biggr)\frac{l^2+k^2}{k^2(l^2-k^2)^2}\Biggr).
\end{equation}

\vspace{7mm}

Using (\ref{101.103}), we write the relation (\ref{101.104}) in 
the following form

\vspace{4mm}

$$
{\sf M}\left\{\left(I_{(01)T,t}
^{*(i_1i_1)}-I_{(01)T,t}^{*(i_1i_1)q}\right)^2\right\}=
{\sf M}\left\{\left(I_{(10)T,t}^{*(i_1i_1)}
-I_{(10)T,t}^{*(i_1i_1)q}\right)^2\right\}=
$$

\vspace{2mm}
$$
=\frac{(T-t)^{4}}{4}\Biggl(\frac{17}{240}-
\frac{1}{3\pi^2}
\sum_{r=1}^q \frac{1}{r^2}-\frac{2}{\pi^4}
\sum_{r=1}^q \frac{1}{r^4} +\Biggr.
$$

\vspace{2mm}
\begin{equation}
\label{daug}
\Biggl.+ \frac{1}{\pi^4}
\Biggl(
\sum_{r=1}^q \frac{1}{r^2}\Biggr)^2-
\frac{1}{\pi^4}
\sum_{\stackrel{k,l=1}{{}_{k\ne l}}}^{q}
\frac{l^2+k^2}{k^2(l^2-k^2)^2}\Biggr).
\end{equation}

\vspace{7mm}

In Tables 1--3, we confirm numerically the formulas 
(\ref{101.100})--(\ref{101.102}), 
(\ref{daug})
for various values $q$. In Tables 1--3,
the number $\varepsilon$  
means the right-hand sides of the mentioned formulas.

The formulas (\ref{101.103}), (\ref{daug1}) appear to be interesting. 
Let us 
confirm numerically their correctness in Tables 4 and 5 
(the number $\varepsilon_q$ is the absolute 
deviation of multiple partial sums with 
the upper limit of summation $q$ for the series (\ref{101.103}), (\ref{daug1})
from the right-hand sides of the formulas (\ref{101.103}), (\ref{daug1});
convergence of multiple series is regarded here 
when $p_1=p_2=q\to\infty$, which is acceptable according to Theorems 1, 2).

\begin{table}
\centering
\caption{Confirmation of the formula (\ref{101.100})}
\label{tab:37}      
\begin{tabular}{p{2.1cm}p{1.7cm}p{1.7cm}p{2.1cm}p{2.3cm}p{2.3cm}p{2.3cm}}
\hline\noalign{\smallskip}
$\varepsilon/(T-t)^3$&0.0459&0.0072&$7.5722\cdot 10^{-4}$
&$7.5973\cdot 10^{-5}$&
$7.5990\cdot 10^{-6}$\\
\noalign{\smallskip}\hline\noalign{\smallskip}
$q$&1&10&100&1000&10000\\
\noalign{\smallskip}\hline\noalign{\smallskip}
\end{tabular}
\end{table}

\begin{table}
\centering
\caption{Confirmation of the formulas (\ref{101.101}), (\ref{101.102})}
\label{tab:38}      
\begin{tabular}{p{2.1cm}p{1.7cm}p{1.7cm}p{2.1cm}p{2.3cm}p{2.3cm}p{2.3cm}}
\hline\noalign{\smallskip}
$4\varepsilon/(T-t)^4$&0.0540&0.0082&$8.4261\cdot 10^{-4}$
&$8.4429\cdot 10^{-5}$&
$8.4435\cdot 10^{-6}$\\
\noalign{\smallskip}\hline\noalign{\smallskip}
$q$&1&10&100&1000&10000\\
\noalign{\smallskip}\hline\noalign{\smallskip}
\end{tabular}
\end{table}

\begin{table}
\centering
\caption{Confirmation of the formula (\ref{daug})}
\label{tab:39}      
\begin{tabular}{p{2.1cm}p{1.7cm}p{1.7cm}p{2.1cm}p{2.3cm}p{2.3cm}p{2.3cm}}
\hline\noalign{\smallskip}
$4\varepsilon/(T-t)^4$&0.0268&0.0034&$3.3955\cdot 10^{-4}$
&$3.3804\cdot 10^{-5}$&
$3.3778\cdot 10^{-6}$\\
\noalign{\smallskip}\hline\noalign{\smallskip}
$q$&1&10&100&1000&10000\\
\noalign{\smallskip}\hline\noalign{\smallskip}
\end{tabular}
\end{table}

\begin{table}
\centering
\caption{Confirmation of the formula (\ref{101.103})}
\label{tab:40}      
\begin{tabular}{p{1.3cm}p{1.8cm}p{1.8cm}p{1.8cm}p{1.8cm}p{2.3cm}p{1.8cm}}
\hline\noalign{\smallskip}
$\varepsilon_q$&2.0294&0.3241&0.0330
&0.0033&
$3.2902\cdot 10^{-4}$\\
\noalign{\smallskip}\hline\noalign{\smallskip}
$q$&1&10&100&1000&10000\\
\noalign{\smallskip}\hline\noalign{\smallskip}
\end{tabular}
\end{table}

\begin{table}
\centering
\caption{Confirmation of the formula (\ref{daug1})}
\label{tab:41}      
\begin{tabular}{p{1.1cm}p{1.5cm}p{1.5cm}p{1.5cm}p{1.5cm}p{1.5cm}p{1.5cm}}
\hline\noalign{\smallskip}
$\varepsilon_q$&10.9585&1.8836&0.1968
&0.0197&
0.0020\\
\noalign{\smallskip}\hline\noalign{\smallskip}
$q$&1&10&100&1000&10000\\
\noalign{\smallskip}\hline\noalign{\smallskip}
\end{tabular}
\end{table}

Using the trigonometric system of functions, let us consider
the approximations of iterated stochastic integrals of the following form

\vspace{-1mm}
$$
{J}_{(\lambda_{1}\ldots \lambda_k)T,t}^{*(i_1\ldots
i_k)}=
{\int\limits_t^{*}}^T\ldots
{\int\limits_t^{*}}^{t_2}
d{\bf w}_{t_{1}}^{(i_1)}\ldots
d{\bf w}_{t_k}^{(i_k)},
$$

\vspace{2mm}
\noindent
where $\lambda_l=1$ if $i_l=1,\ldots,m$ and
$\lambda_l=0$ if $i_l=0;$ $l=1,\ldots,k$
(${\bf w}_{\tau}^{(i)}={\bf f}_{\tau}^{(i)}$
for $i=1,\ldots,m$ and
${\bf w}_{\tau}^{(0)}=\tau$).

It is easy to see that the approximations

\vspace{-1mm}
$$
J_{(\lambda_1\lambda_2)T,t}^{*(i_1 i_2)q},\ \ \
J_{(\lambda_1\lambda_2\lambda_3)T,t}^{*(i_1 i_2 i_3)q}
$$ 

\vspace{2mm}
\noindent
of the stochastic integrals  

\vspace{-1mm}
$$
J_{(\lambda_1\lambda_2)T,t}^{*(i_1 i_2)},\ \ \
J_{(\lambda_1\lambda_2\lambda_3)T,t}^{*(i_1 i_2 i_3)}
$$

\vspace{2mm}
\noindent 
are defined by the right-hand sides of the
formulas (\ref{430}), (\ref{44}),  where it is necessary to take

\vspace{-1mm}
\begin{equation}
\label{123}
\zeta_j^{(i)}=\int\limits_t^T\phi_j(s)d{\bf w}_s^{(i)}
\end{equation}

\vspace{2mm}
\noindent
and $i_1, i_2, i_3=0, 1,\ldots,m$.

Since

\vspace{-1mm}
$$
\int\limits_t^T\phi_j(s)d{\bf w}_s^{(0)}=
\begin{cases}\sqrt{T-t} &\hbox{if}\ j=0\cr\cr
0 &\hbox{if}\ j\ne 0
\end{cases},
$$ 

\vspace{3mm}
\noindent
then
it is easy to get
from (\ref{430}) and (\ref{44}), considering
that in 
these equalities $\zeta_j^{(i)}$ has the form (\ref{123})
and
$i_1, i_2, i_3=0, 1,\ldots,m$,
the following family of formulas

\vspace{3mm}
$$
J_{(10)T,t}^{(i_1 0)q}=
\frac{1}{2}(T-t)^{3/2}\Biggl(
\zeta_0^{(i_1)}+\frac{\sqrt{2}}{\pi}
\Biggl(\sum_{r=1}^{q}\frac{1}{r}\zeta_{2r-1}^{(i_1)}+
\sqrt{\alpha_q}\xi_q^{(i_1)}\Biggr)\Biggr),
$$

\vspace{4mm}
$$
J_{(01)T,t}^{(0 i_2)q}=
\frac{1}{2}(T-t)^{3/2}\Biggl(
\zeta_0^{(i_2)}-\frac{\sqrt{2}}{\pi}
\Biggl(\sum_{r=1}^{q}\frac{1}{r}\zeta_{2r-1}^{(i_2)}+
\sqrt{\alpha_q}\xi_q^{(i_2)}\Biggr)\Biggr),
$$

\vspace{6mm}
$$
J_{(001)T,t}^{(00 i_3)q}=
(T-t)^{5/2}\Biggl(
\frac{1}{6}\zeta_0^{(i_3)}+\frac{1}{2\sqrt{2}\pi^2}\Biggl(
\sum_{r=1}^{q}\frac{1}{r^2}\zeta_{2r}^{(i_3)}+
\sqrt{\beta_q}\mu_q^{(i_3)}\Biggr)-\Biggr.
$$

$$
\Biggl.-
\frac{1}{2\sqrt{2}\pi}\Biggl(\sum_{r=1}^q
\frac{1}{r}\zeta_{2r-1}^{(i_3)}+\sqrt{\alpha_q}\xi_q^{(i_3)}\Biggr)\Biggr),
$$

\vspace{6mm}
$$
J_{(010)T,t}^{(0 i_2 0)q}=
(T-t)^{5/2}\Biggl(
\frac{1}{6}\zeta_0^{(i_2)}-\frac{1}{\sqrt{2}\pi^2}
\Biggl(\sum_{r=1}^{q}\frac{1}{r^2}\zeta_{2r}^{(i_2)}+
\sqrt{\beta_q}\mu_q^{(i_2)}\Biggr)\Biggr),
$$

\vspace{6mm}
$$
J_{(100)T,t}^{(i_1 0 0)q}=
(T-t)^{5/2}\Biggl(
\frac{1}{6}\zeta_0^{(i_1)}+\frac{1}{2\sqrt{2}\pi^2}\Biggl(
\sum_{r=1}^{q}\frac{1}{r^2}\zeta_{2r}^{(i_1)}+
\sqrt{\beta_q}\mu_q^{(i_1)}\Biggr)+\Biggr.
$$

$$
\Biggl.+
\frac{1}{2\sqrt{2}\pi}\Biggl(\sum_{r=1}^q
\frac{1}{r}\zeta_{2r-1}^{(i_1)}+\sqrt{\alpha_q}\xi_q^{(i_1)}\Biggr)\Biggr),
$$

\vspace{9mm}

$$
J_{(011)T,t}
^{*(0 i_2 i_3)q}=(T-t)^{2}\Biggl(\frac{1}{6}
\zeta_{0}^{(i_2)}\zeta_{0}^{(i_3)}-\frac{1}{2\sqrt{2}\pi}
\sqrt{\alpha_q}\xi_q^{(i_3)}\zeta_0^{(i_2)}+\Biggr.
$$

$$
+\frac{1}{2\sqrt{2}\pi^2}\sqrt{\beta_q}\Biggl(
\mu_q^{(i_3)}\zeta_0^{(i_2)}-2\mu_q^{(i_2)}\zeta_0^{(i_3)}\Biggr)+
$$

$$
+\frac{1}{2\sqrt{2}}\sum_{r=1}^{q}
\Biggl(-\frac{1}{\pi r}
\zeta_{2r-1}^{(i_3)}
\zeta_{0}^{(i_2)}+
\frac{1}{\pi^2 r^2}\left(
\zeta_{2r}^{(i_3)}
\zeta_{0}^{(i_2)}-
2\zeta_{2r}^{(i_2)}
\zeta_{0}^{(i_3)}\right)\Biggr)-
$$

$$
-
\frac{1}{2\pi^2}\sum_{\stackrel{r,l=1}{{}_{r\ne l}}}^{q}
\frac{1}{r^2-l^2}\Biggl(
\zeta_{2r}^{(i_2)}
\zeta_{2l}^{(i_3)}+
\frac{l}{r}
\zeta_{2r-1}^{(i_2)}
\zeta_{2l-1}^{(i_3)}
\Biggr)+
$$

$$
+
\sum_{r=1}^{q}
\Biggl(\frac{1}{4\pi r}\left(
\zeta_{2r}^{(i_2)}
\zeta_{2r-1}^{(i_3)}-
\zeta_{2r-1}^{(i_2)}
\zeta_{2r}^{(i_3)}\right)+\Biggr.
$$

\begin{equation}
\label{9000back}
\Biggl.\Biggl.+
\frac{1}{8\pi^2 r^2}\left(
3\zeta_{2r-1}^{(i_2)}
\zeta_{2r-1}^{(i_3)}+
\zeta_{2r}^{(i_3)}
\zeta_{2r}^{(i_2)}\right)\Biggr)\Biggr),
\end{equation}

\vspace{9mm}

$$
J_{(110)T,t}
^{*(i_1 i_2 0)q}=(T-t)^{2}\Biggl(\frac{1}{6}
\zeta_{0}^{(i_1)}\zeta_{0}^{(i_2)}+
\frac{1}{2\sqrt{2}\pi}
\sqrt{\alpha_q}\xi_q^{(i_1)}\zeta_0^{(i_2)}+\Biggr.
$$

$$
+\frac{1}{2\sqrt{2}\pi^2}\sqrt{\beta_q}\Biggl(
\mu_q^{(i_1)}\zeta_0^{(i_2)}-2\mu_q^{(i_2)}\zeta_0^{(i_1)}\Biggr)+
$$

$$
+\frac{1}{2\sqrt{2}}\sum_{r=1}^{q}
\Biggl(\frac{1}{\pi r}
\zeta_{2r-1}^{(i_1)}
\zeta_{0}^{(i_2)}+
\frac{1}{\pi^2 r^2}\left(\zeta_{2r}^{(i_1)}
\zeta_{0}^{(i_2)}
-2\zeta_{2r}^{(i_2)}
\zeta_{0}^{(i_1)}\right)\Biggr)+
$$

$$
+
\frac{1}{2\pi^2}\sum_{\stackrel{r,l=1}{{}_{r\ne l}}}^{q}
\frac{1}{r^2-l^2}\Biggl(\frac{r}{l}
\zeta_{2r-1}^{(i_1)}
\zeta_{2l-1}^{(i_2)}
+\zeta_{2r}^{(i_1)}
\zeta_{2l}^{(i_2)}
\Biggr)+
$$

$$
+
\sum_{r=1}^{q}
\Biggl(\frac{1}{4\pi r}\left(\zeta_{2r-1}^{(i_2)}
\zeta_{2r}^{(i_1)}
-\zeta_{2r-1}^{(i_1)}
\zeta_{2r}^{(i_2)}\right)+\Biggr.
$$

$$
\Biggl.\Biggl.+
\frac{1}{8\pi^2 r^2}\left(
3\zeta_{2r-1}^{(i_1)}
\zeta_{2r-1}^{(i_2)}+
\zeta_{2r}^{(i_1)}
\zeta_{2r}^{(i_2)}\right)\Biggr)\Biggr),
$$

\vspace{9mm}

$$
J_{(101)T,t}
^{*(i_1 0 i_3)q}=(T-t)^{2}\Biggl(\frac{1}{6}
\zeta_{0}^{(i_1)}\zeta_{0}^{(i_3)}
+\frac{1}{2\sqrt{2}\pi}\sqrt{\alpha_q}\left(
\xi_q^{(i_1)}\zeta_0^{(i_3)}-\xi_q^{(i_3)}\zeta_0^{(i_1)}\right)+
\Biggr.
$$

$$
+\frac{1}{2\sqrt{2}\pi^2}\sqrt{\beta_q}\left(
\mu_q^{(i_1)}\zeta_0^{(i_3)}+\mu_q^{(i_3)}\zeta_0^{(i_1)}\right)+
$$

$$
+
\frac{1}{2\sqrt{2}}\sum_{r=1}^{q}
\Biggl(\frac{1}{\pi r}
\left(\zeta_{2r-1}^{(i_1)}
\zeta_{0}^{(i_3)}-\zeta_{2r-1}^{(i_3)}
\zeta_{0}^{(i_1)}\right)+\Biggr.
$$

$$
\Biggl.+
\frac{1}{\pi^2 r^2}\left(
\zeta_{2r}^{(i_1)}
\zeta_{0}^{(i_3)}+
\zeta_{2r}^{(i_3)}
\zeta_{0}^{(i_1)}\right)\Biggr)-
\frac{1}{2\pi^2}\sum_{\stackrel{r,l=1}{{}_{r\ne l}}}^{q}
\frac{1}{rl}
\zeta_{2r-1}^{(i_1)}
\zeta_{2l-1}^{(i_3)}-
$$

$$
\Biggl.-
\sum_{r=1}^q\frac{1}{4\pi^2 r^2}\left(
3\zeta_{2r-1}^{(i_1)}
\zeta_{2r-1}^{(i_3)}+
\zeta_{2r}^{(i_1)}
\zeta_{2r}^{(i_3)}\right)\Biggr).
$$

\vspace{5mm}

\section{Theorems 1--7 from Point
of View of the Wong--Zakai Approximation}

\vspace{5mm}

The iterated Ito stochastic integrals and solutions
of Ito SDEs are complex and important func\-ti\-o\-nals
from the independent components ${\bf f}_{s}^{(i)},$
$i=1,\ldots,m$ of the multidimensional
Wiener process ${\bf f}_{s},$ $s\in[0, T].$
Let ${\bf f}_{s}^{(i)p},$ $p\in\mathbb{N}$ 
be some approximation of
${\bf f}_{s}^{(i)},$
$i=1,\ldots,m$.
Suppose that 
${\bf f}_{s}^{(i)p}$
converges to
${\bf f}_{s}^{(i)},$
$i=1,\ldots,m$ if $p\to\infty$ in some sense and has
differentiable sample trajectories.

A natural question arises: if we replace 
${\bf f}_{s}^{(i)}$
by ${\bf f}_{s}^{(i)p},$
$i=1,\ldots,m$ in the functionals
mentioned above, will the resulting
functionals converge to the original
functionals from the components 
${\bf f}_{s}^{(i)},$
$i=1,\ldots,m$ of the multidimentional
Wiener process ${\bf f}_{s}$?
The answere to this question is negative 
in the general case. However, 
in the pioneering works of Wong E. and Zakai M. \cite{W-Z-1},
\cite{W-Z-2},
it was shown that under the special conditions and 
for some types of approximations 
of the Wiener process the answere is affirmative
with one peculiarity: the convergence takes place 
to the iterated Stratonovich stochastic integrals
and solutions of Stratonovich SDEs and not to iterated 
Ito stochastic integrals and solutions
of Ito SDEs.
The piecewise 
linear approximation 
as well as the regularization by convolution 
\cite{W-Z-1}-\cite{Watanabe} relate the 
mentioned types of approximations
of the Wiener process. The above approximation 
of stochastic integrals and solutions of SDEs 
is often called the Wong--Zakai approximation.

Let ${\bf w}_{\tau},$ $\tau\in[0, T]$ is a random vector with 
an $m+1$ components: ${\bf w}_{\tau}^{(i)}={\bf f}_{\tau}^{(i)}$ 
for $i=1,\ldots,m$ and 
${\bf w}_{\tau}^{(0)}=\tau,$\ 
${\bf f}_{\tau}^{(i)}$ $(i=1,\ldots,m)$
are independent standard Wiener processes.

It is well known that the following representation 
takes place \cite{Lipt}, \cite{7e}

\begin{equation}
\label{um1x}
{\bf w}_{\tau}^{(i)}-{\bf w}_{t}^{(i)}=
\sum_{j=0}^{\infty}\int\limits_t^{\tau}
\phi_j(s)ds\ \zeta_j^{(i)},\ \ \ \zeta_j^{(i)}=
\int\limits_t^T \phi_j(s)d{\bf w}_s^{(i)},
\end{equation}

\vspace{4mm}
\noindent
where $\tau\in[t, T],$ $t\ge 0,$
$\{\phi_j(x)\}_{j=0}^{\infty}$ is an arbitrary complete 
orthonormal system of functions in the space $L_2([t, T]),$ and
$\zeta_j^{(i)}$ are independent standard Gaussian 
random variables for various $i$ or $j.$
Moreover, the series (\ref{um1x}) converges for any $\tau\in [t, T]$
in the mean-square sense.

Let ${\bf w}_{\tau}^{(i)p}-{\bf w}_{t}^{(i)p}$ be 
the mean-square approximation of the process
${\bf w}_{\tau}^{(i)}-{\bf w}_{t}^{(i)},$
which has the following form

\vspace{-3mm}
\begin{equation}
\label{um1xx}
{\bf w}_{\tau}^{(i)p}-{\bf w}_{t}^{(i)p}=
\sum_{j=0}^{p}\int\limits_t^{\tau}
\phi_j(s)ds\ \zeta_j^{(i)}.
\end{equation}

\vspace{3mm}

From (\ref{um1xx}) we obtain

\vspace{-4mm}
\begin{equation}
\label{um1xxx}
d{\bf w}_{\tau}^{(i)p}=
\sum_{j=0}^{p}
\phi_j(\tau)\zeta_j^{(i)} d\tau.
\end{equation}

\vspace{4mm}

Consider the following iterated Riemann--Stieltjes
integral

\begin{equation}
\label{um1xxxx}
\int\limits_t^T
\psi_k(t_k)\ldots \int\limits_t^{t_2}\psi_1(t_1)
d{\bf w}_{t_1}^{(i_1)p_1}\ldots d{\bf w}_{t_k}^{(i_k)p_k},
\end{equation}

\vspace{4mm}
\noindent
where $p_1,\ldots,p_k\in\mathbb{N},$\ \ $i_1,\ldots,i_k=0,1,\ldots,m,$ 

\begin{equation}
\label{um1xxx1}
d{\bf w}_{\tau}^{(i)p}=
\left\{\begin{matrix}
d{\bf f}_{\tau}^{(i)p}\ &\hbox{\rm for}\ \ \ i=1,\ldots,m\cr\cr\cr
d\tau^p\ &\hbox{\rm for}\ \ \ i=0
\end{matrix}
,\right.
\end{equation}

\vspace{4mm}
\noindent
and $d{\bf f}_{\tau}^{(i)p},$ $d\tau^p$ are defined by the relation (\ref{um1xxx}).

Let us substitute (\ref{um1xxx}) into (\ref{um1xxxx})

\begin{equation}
\label{um1xxxx1}
\int\limits_t^T
\psi_k(t_k)\ldots \int\limits_t^{t_2}\psi_1(t_1)
d{\bf w}_{t_1}^{(i_1)p_1}\ldots d{\bf w}_{t_k}^{(i_k)p_k}=
\sum\limits_{j_1=0}^{p_1}\ldots \sum\limits_{j_k=0}^{p_k}
C_{j_k \ldots j_1}\prod\limits_{l=1}^k \zeta_{j_l}^{(i_l)},
\end{equation}

\vspace{4mm}
\noindent
where 
$$
\zeta_j^{(i)}=\int\limits_t^T \phi_j(s)d{\bf w}_s^{(i)}
$$ 

\vspace{2mm}
\noindent
are independent standard Gaussian random variables for various 
$i$ or $j$ (in the case when $i\ne 0$),
${\bf w}_{s}^{(i)}={\bf f}_{s}^{(i)}$ for
$i=1,\ldots,m$ and 
${\bf w}_{s}^{(0)}=s,$

$$
C_{j_k \ldots j_1}=\int\limits_t^T\psi_k(t_k)\phi_{j_k}(t_k)\ldots
\int\limits_t^{t_2}
\psi_1(t_1)\phi_{j_1}(t_1)
dt_1\ldots dt_k
$$

\vspace{4mm}
\noindent
is the Fourier coefficient.

To best of our knowledge \cite{W-Z-1}-\cite{Watanabe}
the approximations of the Wiener process
in the Wong--Zakai approximation must satisfy fairly strong
restrictions
\cite{Watanabe}
(see Definition 7.1, pp.~480--481).
Moreover, approximations of the Wiener process that are
similar to (\ref{um1xx})
were not considered in \cite{W-Z-1}, \cite{W-Z-2}
(also see \cite{Watanabe}, Theorems 7.1, 7.2).
Therefore, the proof of analogs of Theorems 7.1 and 7.2 \cite{Watanabe}
for approximations of the Wiener 
process based on its series expansion (\ref{um1x})
should be carried out separately.
Thus, the mean-square convergence of the right-hand side
of (\ref{um1xxxx1}) to the iterated Stratonovich stochastic integral 
(\ref{str})
does not follow from the results of the papers
\cite{W-Z-1}, \cite{W-Z-2} (also see \cite{Watanabe},
Theorems 7.1, 7.2).

From the other hand, Theorems 1--7 from this 
paper can be considered as the proof of the
Wong--Zakai approximation for the iterated 
Stratonovich stochastic integrals (\ref{str}) of multiplicities 1 to 6
based on the approximation (\ref{um1xx}) of the Wiener process.
At that, the iterated Riemann--Stieltjes integrals (\ref{um1xxxx}) converge
(according to Theorems 1--7)
to the appropriate iterated Stratonovich 
stochastic integrals (\ref{str}). Recall that
$\{\phi_j(x)\}_{j=0}^{\infty}$ (see (\ref{um1x}), (\ref{um1xx}), and
Theorems 3--7)
is a complete 
orthonormal system of Legendre polynomials or 
trigonometric functions 
in the space $L_2([t, T])$.

To illustrate the above reasoning, 
consider two examples for the case $k=2,$
$\psi_1(s),$ $\psi_2(s)\equiv 1;$ $i_1, i_2=1,\ldots,m.$

The first example relates to the piecewise linear approximation
of the multidimensional Wiener process (these approximations 
were considered in \cite{W-Z-1}-\cite{Watanabe}).

Let ${\bf b}_{\Delta}^{(i)}(t),$ $t\in[0, T]$ be the piecewise
linear approximation of the $i$th component ${\bf f}_t^{(i)}$
of the multidimensional standard Wiener process ${\bf f}_t,$
$t\in [0, T]$ with independent components
${\bf f}_t^{(i)},$ $i=1,\ldots,m,$ i.e.

$$
{\bf b}_{\Delta}^{(i)}(t)={\bf f}_{k\Delta}^{(i)}+
\frac{t-k\Delta}{\Delta}\Delta{\bf f}_{k\Delta}^{(i)},
$$

\vspace{4mm}
\noindent
where 

\vspace{-1mm}
$$
\Delta{\bf f}_{k\Delta}^{(i)}={\bf f}_{(k+1)\Delta}^{(i)}-
{\bf f}_{k\Delta}^{(i)},\ \ \
t\in[k\Delta, (k+1)\Delta),\ \ \ k=0, 1,\ldots, N-1.
$$

\vspace{5mm}

Note that w.~p.~1

\begin{equation}
\label{pridum}
\frac{d{\bf b}_{\Delta}^{(i)}}{dt}(t)=
\frac{\Delta{\bf f}_{k\Delta}^{(i)}}{\Delta},\ \ \
t\in[k\Delta, (k+1)\Delta),\ \ \ k=0, 1,\ldots, N-1.
\end{equation}

\vspace{5mm}

Consider the following iterated Riemann--Stieltjes
integral

\vspace{1mm}
$$
\int\limits_0^T
\int\limits_0^{s}
d{\bf b}_{\Delta}^{(i_1)}(\tau)d{\bf b}_{\Delta}^{(i_2)}(s),\ \ \ 
i_1,i_2=1,\ldots,m.
$$

\vspace{5mm}

Using (\ref{pridum}) and additive property of the Riemann--Stieltjes integral, 
we can write w.~p.~1

\vspace{2mm}
$$
\int\limits_0^T
\int\limits_0^{s}
d{\bf b}_{\Delta}^{(i_1)}(\tau)d{\bf b}_{\Delta}^{(i_2)}(s)=
\int\limits_0^T
\int\limits_0^{s}
\frac{d{\bf b}_{\Delta}^{(i_1)}}{d\tau}(\tau)d\tau
\frac{d {\bf b}_{\Delta}^{(i_2)}}{d s}(s)
ds =
$$

\vspace{4mm}
$$
=
\sum\limits_{l=0}^{N-1}\int\limits_{l\Delta}^{(l+1)\Delta}
\left(
\sum\limits_{q=0}^{l-1}\int\limits_{q\Delta}^{(q+1)\Delta}
\frac{\Delta{\bf f}_{q\Delta}^{(i_1)}}{\Delta}d\tau+
\int\limits_{l\Delta}^{s}
\frac{\Delta{\bf f}_{l\Delta}^{(i_1)}}{\Delta}d\tau\right)
\frac{\Delta{\bf f}_{l\Delta}^{(i_2)}}{\Delta}ds=
$$

\vspace{4mm}
$$
=\sum\limits_{l=0}^{N-1}\sum\limits_{q=0}^{l-1}
\Delta{\bf f}_{q\Delta}^{(i_1)}
\Delta{\bf f}_{l\Delta}^{(i_2)}+
\frac{1}{\Delta^2}\sum\limits_{l=0}^{N-1}
\Delta{\bf f}_{l\Delta}^{(i_1)}
\Delta{\bf f}_{l\Delta}^{(i_2)}
\int\limits_{l\Delta}^{(l+1)\Delta}
\int\limits_{l\Delta}^{s}d\tau ds=
$$

\vspace{4mm}

\begin{equation}
\label{oh-ty}
=\sum\limits_{l=0}^{N-1}\sum\limits_{q=0}^{l-1}
\Delta{\bf f}_{q\Delta}^{(i_1)}
\Delta{\bf f}_{l\Delta}^{(i_2)}+
\frac{1}{2}\sum\limits_{l=0}^{N-1}
\Delta{\bf f}_{l\Delta}^{(i_1)}
\Delta{\bf f}_{l\Delta}^{(i_2)}.
\end{equation}

\vspace{6mm}

Using (\ref{oh-ty}), it 
is not difficult to show 
that

\vspace{1mm}
$$
\hbox{\vtop{\offinterlineskip\halign{
\hfil#\hfil\cr
{\rm l.i.m.}\cr
$\stackrel{}{{}_{N\to \infty}}$\cr
}} }
\int\limits_0^T
\int\limits_0^{s}
d{\bf b}_{\Delta}^{(i_1)}(\tau)d{\bf b}_{\Delta}^{(i_2)}(s)=
\int\limits_0^T
\int\limits_0^{s}
d{\bf f}_{\tau}^{(i_1)}d{\bf f}_{s}^{(i_2)}+
\frac{1}{2}{\bf 1}_{\{i_1=i_2\}}\int\limits_0^T ds=
$$

\vspace{4mm}
\begin{equation}
\label{uh-111}
=
\int\limits_0^{*T}
\int\limits_0^{*s}
d{\bf f}_{\tau}^{(i_1)}d{\bf f}_{s}^{(i_2)},
\end{equation}

\vspace{6mm}
\noindent
where $\Delta\to 0$ if $N\to\infty$ ($N\Delta=T$).

Obviously, (\ref{uh-111}) agrees with Theorem 7.1 (see \cite{Watanabe},
p.~486).

The next example relates to the approximation (\ref{um1xx})
of the Wiener process based on its series expansion
(\ref{um1x}) for $t=0$, where
$\{\phi_j(x)\}_{j=0}^{\infty}$ 
is a complete 
orthonormal system of Legendre polynomials or 
trigonometric functions 
in the space $L_2([0, T])$.

Consider the following iterated Riemann--Stieltjes
integral

\begin{equation}
\label{abcd1}
\int\limits_0^T
\int\limits_0^{s}
d{\bf f}_{\tau}^{(i_1)p}d{\bf f}_{s}^{(i_2)p},\ \ \ 
i_1,i_2=1,\ldots,m,
\end{equation}

\vspace{4mm}
\noindent
where $d{\bf f}_{\tau}^{(i)p}$ is defined by the
relation
(\ref{um1xxx}).

Let us substitute (\ref{um1xxx}) into (\ref{abcd1})

\begin{equation}
\label{set18}
\int\limits_0^T
\int\limits_0^{s}
d{\bf f}_{\tau}^{(i_1)p}d{\bf f}_{s}^{(i_2)p}=
\sum\limits_{j_1,j_2=0}^p
C_{j_2 j_1} \zeta_{j_1}^{(i_1)}\zeta_{j_2}^{(i_2)},
\end{equation}

\vspace{4mm}
\noindent
where 

\vspace{-1mm}
$$
C_{j_2 j_1}=
\int\limits_0^T \phi_{j_2}(s)\int\limits_0^s
\phi_{j_1}(\tau)d\tau ds
$$

\vspace{4mm}
\noindent
is the Fourier coefficient; another notations 
are the same as in (\ref{um1xxxx1}).

As we noted above, approximations of the Wiener process that are
similar to (\ref{um1xx})
were not considered in \cite{W-Z-1}, \cite{W-Z-2}
(also see Theorems 7.1, 7.2 in \cite{Watanabe}).
Furthermore, the extension of the results of Theorems 7.1 and 7.2
\cite{Watanabe} to the case under consideration is
not obvious.

On the other hand, we can apply the theory built in Chapters 1 and 2
of the monographs \cite{2018a}-\cite{2018axx}. More precisely, 
using 
Theorem 3, we obtain from (\ref{set18}) the desired result

\vspace{1mm}
$$
\hbox{\vtop{\offinterlineskip\halign{
\hfil#\hfil\cr
{\rm l.i.m.}\cr
$\stackrel{}{{}_{p\to \infty}}$\cr
}} }
\int\limits_0^T
\int\limits_0^{s}
d{\bf f}_{\tau}^{(i_1)p}d{\bf f}_{s}^{(i_2)p}=
\hbox{\vtop{\offinterlineskip\halign{
\hfil#\hfil\cr
{\rm l.i.m.}\cr
$\stackrel{}{{}_{p\to \infty}}$\cr
}} }
\sum\limits_{j_1,j_2=0}^p
C_{j_2 j_1} \zeta_{j_1}^{(i_1)}\zeta_{j_2}^{(i_2)}=
$$

\vspace{3mm}
\begin{equation}
\label{umen-bl}
=
\int\limits_0^{*T}
\int\limits_0^{*s}
d{\bf f}_{\tau}^{(i_1)}d{\bf f}_{s}^{(i_2)}.
\end{equation}

\vspace{7mm}

From the other hand, by Theorems 1, 2
(see (\ref{a2})) for the case
$k=2$ we obtain from (\ref{set18}) the following relation

$$
\hbox{\vtop{\offinterlineskip\halign{
\hfil#\hfil\cr
{\rm l.i.m.}\cr
$\stackrel{}{{}_{p\to \infty}}$\cr
}} }
\int\limits_0^T
\int\limits_0^{s}
d{\bf f}_{\tau}^{(i_1)p}d{\bf f}_{s}^{(i_2)p}=
\hbox{\vtop{\offinterlineskip\halign{
\hfil#\hfil\cr
{\rm l.i.m.}\cr
$\stackrel{}{{}_{p\to \infty}}$\cr
}} }
\sum\limits_{j_1,j_2=0}^p
C_{j_2 j_1} \zeta_{j_1}^{(i_1)}\zeta_{j_2}^{(i_2)}=
$$

\vspace{5mm}
$$
=
\hbox{\vtop{\offinterlineskip\halign{
\hfil#\hfil\cr
{\rm l.i.m.}\cr
$\stackrel{}{{}_{p\to \infty}}$\cr
}} }
\sum\limits_{j_1,j_2=0}^p
C_{j_2 j_1} \biggl(\zeta_{j_1}^{(i_1)}\zeta_{j_2}^{(i_2)}-
{\bf 1}_{\{i_1=i_2\}}{\bf 1}_{\{j_1=j_2\}}\biggr)+
{\bf 1}_{\{i_1=i_2\}}\sum\limits_{j_1=0}^{\infty}
C_{j_1 j_1}=
$$

\vspace{3mm}

\begin{equation}
\label{umen-blx}
=
\int\limits_0^T
\int\limits_0^{s}
d{\bf f}_{\tau}^{(i_1)}d{\bf f}_{s}^{(i_2)}+
{\bf 1}_{\{i_1=i_2\}}\sum\limits_{j_1=0}^{\infty}
C_{j_1 j_1}.
\end{equation}

\vspace{7mm}

Since

\vspace{1mm}
$$
\sum\limits_{j_1=0}^{\infty}
C_{j_1 j_1}=\frac{1}{2}\sum\limits_{j_1=0}^{\infty}
\left(\int\limits_0^T \phi_j(\tau)d\tau\right)^2
=\frac{1}{2}
\left(\int\limits_0^T \phi_0(\tau)d\tau\right)^2=\frac{1}{2}
\int\limits_0^T ds,
$$

\vspace{7mm}
\noindent
then from (\ref{umen-blx}) we obtain (\ref{umen-bl}).

\end{document}